\documentclass{msripub}
\usepackage{book37}
\usepackage{msribib}
\usepackage{rot}

\def\Maple/{Maple}

\def\figdir/{propp/}
\def\yij{\Young{\mathstrut{\scriptstyle i}&\mathstrut{\scriptstyle j}\cr}}
\def\yIJ{\Young{\mathstrut{\scriptstyle i}\cr\mathstrut{\scriptstyle j}\cr}}

\newtheorem{problem}{Problem}
\newtheorem[{\setprovedbox}{\it\specialdigits}]{theorem}{Theorem}
\unnumbered{Theorem}
\newenvironment{progress}{\proof[Progress]}{\setbox0\hbox{\box\provedbox} \endproof}

\newbox\makex \newbox\makev \newbox\makea

\setbox\makex\vbox{\vskip.6667pt
\hbox{\hskip .6667pt\vrule height 4.4402pt width 4.4402pt\hskip .6667pt}
\vskip.6667pt}

\newdimen\unit \unit=5.7735 pt  % baselineskip / sqrt3

\expandafter\setbox\expandafter\makev\psfig{file=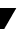}
\expandafter\setbox\expandafter\makea\psfig{file=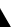}

\def\docats{\catcode`\V=13 \catcode`\A=13 \catcode`\X=13%
        \catcode`\.=13\obeylines \obeyspaces}
\def\undocats{\catcode`\V=11\catcode`\A=11\catcode`\X=11%
        \catcode`\.=12\catcode`\^^M=5\catcode`\ =10}

\docats
\defA{\copy\makea}%
\def.{\hskip\unit}%
\defV{\copy\makev}%
\defX{\copy\makex}%
\def {.}%
\def\newdefs{%
\defV{\hbox{\tt\char`\V}}%
\defA{\hbox{\tt\char`\A}}%
\def.{\hskip5,25pt}%
}
\undocats%

\def\cboard{\vcenter\brd}

\def\brda{\bgroup \let\par=\cr \docats \newdefs
        \baselineskip 10pt \afterassignment\begbrd \def\theboard}

\def\brd{\bgroup \let\par=\cr \docats
        \offinterlineskip \afterassignment\begbrd \def\theboard}
\def\begbrd{\halign{&##\cr\theboard\crcr}\egroup}

\title
[Enumeration of Matchings: Problems and Progress]
{Enumeration of Matchings:\\ Problems and Progress}

\author{James Propp}
\address{Department of Mathematics\\
University of Wisconsin\\
Madison, WI 53706\\
United States}
\email{propp@math.wisc.edu}

\begin{document}

\hrule height 0pt
\vfil
\begin{centering}
\large
{\bf Enumeration of Matchings:\\ Problems and Progress}
\vskip12pt
James Propp
\vskip12pt
\normalsize

{\it New Perspectives in Geometric Combinatorics}

Edited by
L. Billera,
A. Bj\"orner,
C. Greene,
R. Simeon, and R. Stanley

Mathematical Sciences Research Institute Publications {\bf 38}

Cambridge University Press, 1999

Pages 255--291

\end{centering}
\vfil
\eject

\setcounter{page}{254}
\cleardoublepage

\dedicatory{Dedicated to the memory of David Klarner (1940--1999)}

\begin{abstract}
This document is built around a list
of thirty-two problems in enumeration of matchings,
the first twenty of which were presented in a lecture at MSRI
in the fall of 1996.
I begin with a capsule history
of the topic of enumeration of matchings.
The twenty original problems, with commentary,
comprise the bulk of the article.
I give an account of the progress that has been made
on these problems as of this writing,
and include pointers to both the printed and on-line literature;
roughly half of the original twenty problems
were solved by participants in the MSRI Workshop on Combinatorics,
their students, and others, between 1996 and 1999.
The article concludes with a dozen new open problems.
\end{abstract}

\maketitle

\section{Introduction}

How many perfect matchings does a given graph $G$ have?
That is, in how many ways can one choose a subset of the edges of $G$
so that each vertex of $G$ belongs to one and only one chosen edge?
(See Figure 1(a) for an example of a perfect matching of a graph.)
For general graphs $G$,
it is computationally hard to obtain the answer \cite{PrV},
and even when we have the answer, it is not so clear
that we are any the wiser for knowing this number.
However, for
many infinite families of special graphs the number of perfect matchings
is given by compellingly simple formulas.
Over the past ten years a great many families of this kind
have been discovered,
and while there is no single unified result that encompasses all of them,
many of these families resemble one another,
both in terms of the form of the results
and in terms of the methods that have been useful in proving them.

\begin{figure}
\centerline{\psfig{file=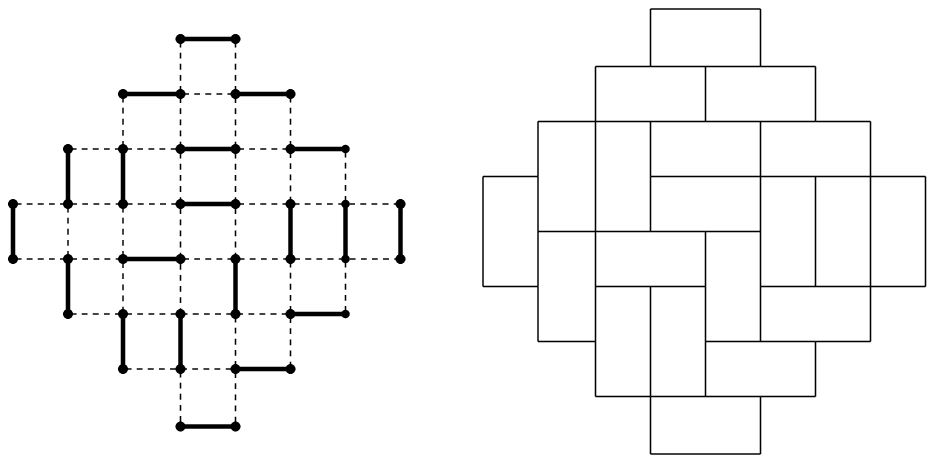,scale=110}}
\caption{The Aztec diamond of order 4.}
\label{azgraph}
\end{figure}

The deeper significance of these formulas is not clear.
Some of them are related to results in representation theory
or the theory of symmetric functions,
but others seem to be self-contained combinatorial puzzles.
Much of the motivation for this branch of research
lies in the fact that we are still unable to predict ahead of time
which enumerative problems lead to beautiful formulas and which do not;
each new positive result seems like an undeserved windfall.

Hereafter, I will use the term ``matching'' to signify ``perfect matching''.
(See the book of Lov\'asz and Plummer \citeyear{PrLP}
for general background on the theory of matchings.)

As far as I have been able to determine, problems involving enumeration
of matchings were first examined by chemists and physicists in the 1930s,
for two different (and unrelated) purposes: the study of aromatic
hydrocarbons and the attempt to create a theory of the liquid state.

Shortly after the advent of quantum chemistry, chemists turned their
attentions to molecules like benzene composed of carbon rings with
attached hydrogen atoms.  For these researchers, matchings of a graph
corresponded to ``Kekul\'e
structures'', i.e., ways of assigning single and double bonds in the
associated hydrocarbon (with carbon atoms at the vertices and tacit
hydrogen atoms attached to carbon atoms with only two neighboring
carbon atoms).  See for example the article of Gordon and Davison
\citeyear{PrGD},
whose use of nonintersecting lattice paths anticipates certain
later work \cite{PrGV,PrSac1,PrJS}.  There
are strong connections between combinatorics and chemistry for such
molecules; for instance, those edges which are present in comparatively
few of the matchings of a graph turn out to correspond to the bonds that
are least stable, and the more matchings a polyhex graph possesses the
more stable is the corresponding benzenoid molecule.  Since hexagonal
rings are so predominant in the structure of hydrocarbons, chemists gave
most of their attention to counting matchings of subgraphs of the infinite
honeycomb grid.

At approximately the same time, scientists were trying to understand
the behavior of liquids.  As an extension of a more basic model for
liquids containing only molecules of one type, Fowler and Rushbrooke
\citeyear{PrFR} devised a lattice-based model for liquids containing two
types of molecules, one large and one small.  In the case where the
large molecule was roughly twice the size of the small molecule, it
made sense to model the small molecules as occupying sites of a
three-dimensional grid and the large molecules as occupying pairs of
adjacent sites.  In modern parlance, this is a monomer-dimer model.
In later years, the two-dimensional version of the model was found
to have applicability to the study of molecules adsorbed on films;
if the adsorption sites are assumed to form a lattice, and an adsorbed
molecule is assumed to occupy two such sites, then one can imagine
fictitious molecules that occupy all the unoccupied sites (one each).

Major progress was made when Temperley and Fisher \citeyear{PrTF}
and Kasteleyn \citeyear{PrKa1} independently found ways to count pure
dimer configurations on subgraphs of the infinite square grid,
with no monomers present.  Although the physical significance of
this special case was (and remains) unclear, this result, along with
Onsager's earlier exact solution of the two-dimensional Ising model
\cite{PrO}, paved the way for other advances such as Lieb's exact solution
of the six-vertex model \cite{PrL}, culminating in a new field at the
intersection of physics and mathematics: exactly solved
statistical mechanics models in two-dimensional lattices.  (Intriguingly,
virtually none of the three- and higher-dimensional analogues of these
models have succumbed to researchers' efforts at obtaining exact solutions.)
For background on lattice models in statistical mechanics, see
the book by Baxter \citeyear{PrB}.

An infinite two-dimensional grid has many finite subgraphs; in choosing
which ones to study, physicists were guided by
the idea that the shape of boundary should be chosen so as to minimize
the effect of the boundary\emdash that is, to maximize the number of
configurations, at least in the asymptotic sense.  For example,
Kasteleyn, in his study of the dimer model on the square grid,
counted the matchings of the $m$-by-$n$ rectangle (see the double-product
formula at the beginning of Section~5) and of the $m$-by-$n$ rectangular torus,
and showed that the two numbers grow at the
same rate as $m,n$ go to infinity, namely $C^{mn}$ for a known constant $C$.
(Analytically, $C$ is $e^{G/\pi}$, where $G$ is Catalan's constant
$1 - \frac{1}{9} + \frac{1}{25} - \frac{1}{49} + \frac{1}{81} - \cdots{}$;
numerically, $C$ is approximately $1.34$.)

Kasteleyn \citeyear{PrKa1} wrote: ``The effect of boundary
conditions is, however, not
entirely trivial and will be discussed in more detail in a subsequent
paper.''  (See the article of Cohn, Kenyon and Propp \cite{PrCKP}
for a rigorous mathematical treatment of boundary conditions.)
Kasteleyn never wrote such a followup paper, but other physicists did
give some attention to the issue of boundary shape, most notably
Grensing, Carlsen and Zapp \cite{PrGCZ}. These authors
considered a one-parameter family of graphs of the kind
shown in Figure 1(a), and they asserted that every graph in this family
has $2^{N/4}$ matchings, where $N$
is the number of vertices.  They did not give a proof, nor did they indicate
whether they had one.  The result was rediscovered in the late 1980s by
Elkies, Kuperberg, Larsen, and Propp \cite{PrEKLP}, who gave four proofs of
the formula.  This article led to a great deal of work among enumerative
combinatorialists, who refer to graphs like the one shown in Figure 1
as ``Aztec diamond graphs'', or sometimes just Aztec diamonds for short.
(It should be noted that Elkies et al.\ \citeyear{PrEKLP} used the term ``Aztec
diamond'' to denote regions like the one shown in Figure 1(b).  The two
sorts of Aztec diamonds are dual to one another; matchings of Aztec diamond
graphs correspond to domino tilings of Aztec diamond regions.)

At about the same time, it became clear that there had been earlier work
within the combinatorial community that was pertinent to the study of
matchings, though its relevance had not hitherto been recognized.  For
instance, Mills, Robbins and Rumsey \cite{PrMRR}, in their work on alternating
sign matrices, had counted pairs of ``compatible'' ASMs of consecutive size;
these can be put into one-to-one correspondence with matchings
of an associated Aztec diamond graph \cite{PrEKLP}.

Looking into earlier mathematical literature,
one can even see intimations of enumerative matching theory in
the work of MacMahon \citeyear{PrM}, who nearly a century ago found a formula
for the number of plane partitions whose solid Young diagram fits inside
an $a$-by-$b$-by-$c$ box, as will be discussed in Section~2.
(See the book by Andrews
\citeyear{PrA} and the article by Stanley \citeyear{PrSt1} for
background on plane partitions.)
Such a Young diagram is nothing more than an assemblage of cubes,
and it has long been known in the extra-mathematical world
that such assemblages, viewed from a distant point, looks like tilings
(consider Islamic art, for instance).  Thus it was natural for
mathematicians to interpret MacMahon's theorem on plane partitions
as a result about tilings of a hexagon by rhombuses.
This insight may have occurred to a number of people independently;
the earliest chain of oral communication that I have followed
leads back to Klarner (who did not publish his observation
but relayed it to Stanley in the 1970s),
and the earliest published statement I have found is in a paper
by David and Tomei \citeyear{PrDT}.{\looseness=-1\par}

In any case, each of the Young diagrams enumerated by MacMahon
corresponds to a tiling of a hexagon by rhombuses,
where the hexagon is semiregular
(its opposite sides are parallel and of equal length,
with all internal angles equal to 120 degrees)
and has side-lengths $a,b,c,a,b,c$,
and where the rhombuses have all side-lengths equal to 1.
These tilings in turn correspond to matchings of the ``honeycomb''
graph that is dual to the dissection of the hexagon into unit equilateral
triangles; see Figure 2, which shows a matching of the honeycomb graph
and the associated tiling of a hexagon.
Kuperberg \citeyear{PrKu1} was the first to exploit the connection
between plane partitions and the dimer model.
(Interestingly, some of the same graphs
that Kuperberg studied had been investigated independently by chemists
in their study of benzenoids hydrocarbons; Cyvin and Gutman \citeyear{PrCG}
give a survey of this work.)

\begin{figure}
\centerline{\psfig{file=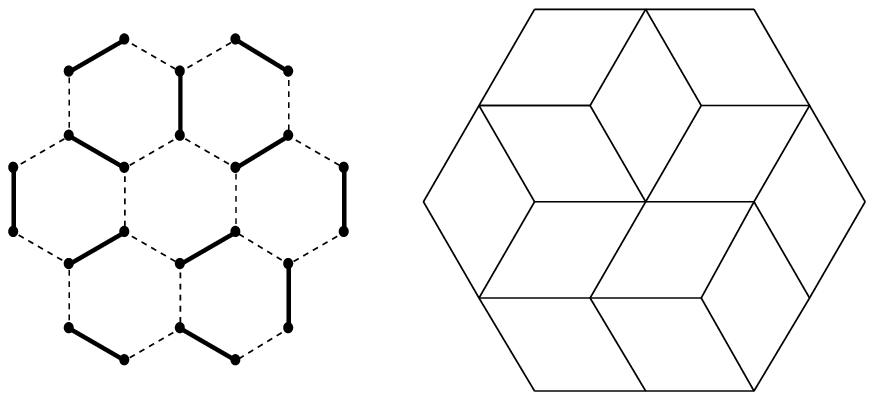,scale=110}}
\caption{A matching and its associated tiling.}
\label{hexgraph}
\end{figure}

Similarly, variants of MacMahon's problem in which the plane
partition is subjected to various symmetry constraints
(considered by Macdonald, Stanley, and others
[\citeNP{PrSt3}; \citeyearNP{PrSt4}])
correspond to the problem of enumerating matchings possessing
corresponding kinds of symmetry.  Kuperberg \citeyear{PrKu1} used
this correspondence in solving one of Stanley's open problems, and this
created further interest in matchings among combinatorialists.

One of Kuperberg's chief tools was an old result of Kasteleyn, which
showed that for any planar graph $G$, the number of matchings
of $G$ is equal to the Pfaffian of a certain matrix of zeros and ones
associated with $G$.  A special case of this result, enunciated by
Percus \citeyear{PrPe},
can be used when $G$ is bipartite; in this case, one can use a determinant
instead of a Pfaffian.  Percus' determinant is a modified version of
the bipartite adjacency matrix of the graph, in which rows correspond to
``white'' vertices and columns correspond to ``black'' vertices (under
a coloring scheme whereby white vertices have only black neighbors and
vice versa); the $(i,j)$-th entry is $\pm 1$ if the $i$-th white vertex and
$j$-th black vertex are adjacent, and 0 otherwise.  For more details on
how the signs of the entries are chosen, see the
expositions of Kasteleyn \citeyear{PrKa3} and Percus \citeyear{PrPe}.

Percus' theorem, incorporated into computer software, makes it easy to
count the matchings of many planar graphs and look for patterns in
the numbers that arise.  Two such programs are {\tt vaxmaple}, written
by Greg Kuperberg, David Wilson and myself, and {\tt vaxmacs}, written
by David Wilson.  Most of the patterns described below were
discovered with the aid of this software, which is available from
\url{http://math.wisc.edu/~propp/software.html}.  Both programs
treat subgraphs of the infinite square grid; this might seem restrictive,
but it turns out that counting the matchings of an arbitrary bipartite
planar graph can be fit into this framework, with a bit of tweaking.
The mathematically interesting part of each program is the routine for
choosing the signs of the nonzero entries.  There are many choices
that would work, but Wilson's sign-rule is far and away the simplest:
If an edge is horizontal, we give it weight $+1$, and if an edge is
vertical, joining a vertex in one row to a vertex in the row below it,
we give the edge weight $(-1)^k$, where $k$ is the number of vertices in
the upper row to the left of the vertical edge.

The main difference between {\tt vaxmaple} and {\tt vaxmacs} is that
the former creates \Maple/ code which, if sent to \Maple/,
results in \Maple/ printing out the number of matchings of the graph;
{\tt vaxmacs}, on the other hand, is a customized {\sc Emacs} environment
that fully integrates text-editing operations (used for defining the
graph one wishes to study) with the mathematical operations of interest.
Both programs represent bipartite planar graphs in ``VAX-format'',
where V's, A's, X's, and other letters denote vertices.  (An example
of VAX-format can be found on page~\pageref{vaxexample}; for a detailed
explanation see \url{http://math.wisc.edu/~propp/vaxmaple.doc}.)

Quite recently, the study of matchings of nonbipartite graphs
has been expedited by the programs {\tt graph} and {\tt planemaple},
created by Matt Blum and Ben Wieland, respectively.
These programs make it easy to define a planar graph by pointing and clicking,
after which one can count its matchings using an efficient implementation
of Kasteleyn's Pfaffian method.  This makes it easy to try out new ideas
and look for patterns, outside of the better-explored bipartite
case.{\looseness=1\par} 

Interested readers with access to the World Wide Web
can obtain copies of all of these programs via
\url{http://math.wisc.edu/~propp/software.html}.

Most of the formulas that have been discovered
express the number of matchings of a graph
as a product of many comparatively small factors.
Even before one has conjectured (let alone proved) such a formula,
one can frequently infer its existence
from the fact that the number of matchings
has only small prime factors.
Numbers that are large compared to their largest prime factor
are sometimes called ``smooth'' or ``round'';
the latter term will be used here.
The definition of roundness is not precise,
since it is not intended for use as a technical term.
Its vagueness is intended to capture the uncertainties and the suspense
of formula-hunting,
and the debatable issue of whether the occurrence of a single
larger-than-expected prime factor
rules out the existence of a product formula.
(For an example of a number whose roundness lies in this gray area,
see the table of numbers given in Problem 8.)
It is worth noting that Kuperberg \citeyear[Section VII-A]{PrKu2}
has shown that rigorous proofs of roundness need not always yield explicit
product formulas.

Christian Krattenthaler has written a Mathematica
program called RATE that greatly expedites the process of guessing patterns
in experimental data on enumeration of matchings; see
\url{http://radon.mat.univie.ac.at/People/kratt/rate/rate.html}.

A great source of the appeal of research on enumeration of matchings is the
ease with which undergraduate research assistants can participate in the hunt
for formulas and proofs; many members of the M.I.T.\ Tilings Research Group
(composed mostly of undergraduates like Blum and Wieland) played a role in
the developments that led to the writing of this article.  Enumeration of
matchings has turned out to be a rich avenue of combinatorial inquiry,
and many more beautiful patterns undoubtedly await discovery.

Updates on the status of these problems can be found on the Web at
\url{http://math.wisc.edu/~propp/update.ps.gz}.

\section{Lozenges}

We begin with problems related to lozenge tilings of hexagons.
A \textit{lozenge} is a rhombus of side-length 1
whose internal angles measure 60 and 120 degrees;
all the hexagons we will consider will tacitly have
integer side-lengths and internal angles of 120 degrees.
Every such hexagon $H$ can be dissected into unit equilateral triangles
in a unique way,
and one can use this dissection to define a graph $G$
whose vertices correspond to the triangles
and whose edges correspond to pairs of triangles that share an edge;
this is the ``finite honeycomb graph''
dual to the dissection.
It is easy to see
that the tilings of $H$ by lozenges are in one-to-one correspondence
with the matchings of $G$.

The $a,b,c$ semiregular hexagon
is the hexagon whose side lengths are, in cyclical order, $a,b,c,a,b,c$.
Lozenge tilings of this region
are in correspondence with
plane partitions with
at most $a$ rows,
at most $b$ columns,
and no part exceeding $c$.
Such hexagons are represented in VAX-format by diagrams like
$$
\label{vaxexample}
\vcenter\brda{%
   AVAVAVAVA
  AVAVAVAVAVA
 AVAVAVAVAVAVA
AVAVAVAVAVAVAVA
VAVAVAVAVAVAVAV
 VAVAVAVAVAVAV
  VAVAVAVAVAV
   VAVAVAVAV
}
$$
where A's and V's represent upward-pointing and downward-pointing triangles,
respectively.  In this article we will use triangles instead:
$$\cboard{
   AVAVAVAVA
  AVAVAVAVAVA
 AVAVAVAVAVAVA
AVAVAVAVAVAVAVA
VAVAVAVAVAVAVAV
 VAVAVAVAVAVAV
  VAVAVAVAVAV
   VAVAVAVAV
}
$$

MacMahon \citeyear{PrM} showed that the number of such plane partitions is
$$\prod_{i=0}^{a-1} \prod_{j=0}^{b-1} \prod_{k=0}^{c-1}
\frac{i+j+k+2}{i+j+k+1}\,.$$
(This form of MacMahon's formula is due to Macdonald; a short,
self-contained proof is given by Cohn et al. \citeyear[Section~2]{PrCLP}.)

\begin{problem}
\label{prob1}
Show that in the $2n-1$, $2n$, $2n-1$ semiregular hexagon,
the central location (consisting of the two innermost triangles) is
covered by a lozenge in exactly one-third of the tilings.

(Equivalently: Show that if one chooses a random matching of
the dual graph, the probability that the central edge is contained
in the matching is exactly $\frac13$.)
\end{problem}

\begin{progress}
Two independent and very different solutions of this problem have been
found; one by Mihai Ciucu and Christian Krattenthaler
and the other by Harald Helfgott and Ira Gessel.
Ciucu and Krattenthaler \citeyear{PrCK}
compute more generally the number of rhombus tilings of a hexagon
with sides $a,a,b,a,a,b$ that contain the central unit rhombus, where $a$
and $b$ must have opposite parity (the special case $a=2n-1$, $b=2n$ solves
Problem 1).  The same generalization was obtained (in a different but
equivalent form) by Helfgott and Gessel \citeyear{PrHG}, using a completely
different method.
One might still try to look for a proof whose simplicity
is comparable to that of the answer ``one-third''.
Also worthy of note is the paper of Fulmek and Krattenthaler \citeyear{PrFK1},
which generalizes the result of Ciucu and Krattenthaler \citeyear{PrCK}.
\end{progress}

\hskip\parindent
The hexagon of side-lengths
$n$, $n+1$, $n$, $n+1$, $n$, $n+1$ cannot be tiled by
lozenges at all, for in the dissection into unit triangles, the number of
upward-pointing triangles differs from the number of downward-pointing
triangles.  However, if one removes the central triangle, one gets a
region that can be tiled, and the sort of numbers one gets for small
values of $n$ are striking.  Here they are, in factored form:
$$
\medmuskip 1mu
\displaylines{
2\cr
2\cdot3^3\cr
2^5\cdot3^3\cdot5\cr
2^5\cdot5^7\cr
2^2\cdot5^7\cdot7^5\cr
2^8\cdot3^3\cdot5\cdot7^{11}\cr
2^{13}\cdot3^9\cdot7^{11}\cdot11\cr
2^{13}\cdot3^{18}\cdot7^5\cdot11^7\cr
2^8\cdot3^{18}\cdot11^{13}\cdot13^5\cr
2^2\cdot3^9\cdot11^{19}\cdot13^{11}\cr
2^{10}\cdot3^3\cdot11^{19}\cdot13^{17}\cdot17\cr
2^{16}\cdot11^{13}\cdot13^{23}\cdot17^7\cr
}
$$
These are similar to the numbers one gets from counting lozenge tilings
of an $n,n,n,n,n,n$ hexagon, in that the largest prime factor seems to
be bounded by a linear function of $n$.

\begin{problem}
\label{prob2}
Enumerate the lozenge tilings of the region obtained
from the $n$, $n+1 $, $n $, $n+1 $, $n $, $n+1$ hexagon by removing
the central triangle.
\end{problem}

\begin{progress}
Mihai Ciucu has solved the more general problem of counting the rhombus
tilings of an $(a$, $b+1 $, $b $, $a+1 $, $b $, $b+1)$-hexagon with the central triangle removed
\cite{PrCi2}.
Ira Gessel proved this result independently using the nonintersecting
lattice-paths method \cite{PrHG}.
Soichi Okada and Christian Krattenthaler have solved the even more general
problem of counting the rhombus tilings of an $(a $, $b+1 $, $c $, $a+1 $, $b $, $c+1)$-hexagon
with the central triangle removed \cite{PrOK}.
\end{progress}

\hskip\parindent
One can also take a $2n $, $2n+3 $, $2n $, $2n+3 $, $2n $, $2n+3$ hexagon
and make it lozenge-tilable by removing a triangle from the
middle of each of its three long sides, as shown:
$$\cboard{
          AVAVAVAVAVAVAVAVA
         AVAVAVAVAVAVAVAVAVA
        AVAVAVAVAVAVAVAVAVAVA
       AVAVAVAVAVAVAVAVAVAVAVA
      AVAVAVAVAVAVAVAVAVAVAVAVA
      VAVAVAVAVAVAVAVAVAVAVAVAV
    AVAVAVAVAVAVAVAVAVAVAVAVAVAVA
   AVAVAVAVAVAVAVAVAVAVAVAVAVAVAVA
  AVAVAVAVAVAVAVAVAVAVAVAVAVAVAVAVA
 AVAVAVAVAVAVAVAVAVAVAVAVAVAVAVAVAVA
AVAVAVAVAVAVAVAVAVAVAVAVAVAVAVAVAVAVA
VAVAVAVAVAVAVAVAVAVAVAVAVAVAVAVAVAVAV
 VAVAVAVAVAVAVAVAVAVAVAVAVAVAVAVAVAV
  VAVAVAVAVAVAVAVAVAVAVAVAVAVAVAVAV
   VAVAVAVAVAVAVAVAVAVAVAVAVAVAVAV
    VAVAVAVAVAVAVAVAVAVAVAVAVAVAV
     VAVAVAVAVAVAVAVAVAVAVAVAVAV
      VAVAVAVAVAVAVAVAVAVAVAVAV
       VAVAVAVAVAV VAVAVAVAVAV
}$$
Here one obtains an equally tantalizing sequence of factorizations:
$$
\medmuskip 1mu
\displaylines{
1\cr
2^7\cdot 7^2\cr
2^2\cdot 7^4\cdot 11^4\cdot 13^2\cr
2^{10}\cdot 3^3\cdot 5^8\cdot 13^2\cdot 17^4 \cdot 19^2\cr
2^2\cdot 5^2\cdot 7^2\cdot 11^3\cdot 13^4\cdot 17^4\cdot 19^8\cdot 23^4
}
$$

\begin{problem}
\label{prob3}
Enumerate the lozenge tilings of the region obtained
from the $2n$, $2n+3 $, $2n $, $2n+3 $, $2n $, $2n+3$ hexagon by
removing a triangle from
the middle of each of its long sides.
\end{problem}

\begin{progress}
Theresia Eisenk\"olbl solved this problem.  What she does in fact is
to compute the number of all rhombus tilings of a hexagon with sides
$a$, $b+3$, $c$, $a+3$, $b$, $c+3$, where an arbitrary triangle is
removed from each of the ``long'' sides of the hexagon (not
necessarily the triangle in the middle).  For the proof of her formula
\cite{PrE1} she uses nonintersecting lattice paths, determinants, and
the Jacobi determinant formula \cite{PrT}. However, I still know of no
conceptual explanation for why these numbers are so close (in the
multiplicative sense) to being perfect squares.
\end{progress}

\hskip\parindent
We now return to ordinary $a,b,c$ semiregular hexagons.
When $a=b=c$, there are not two but six central triangles.
There are two geometrically distinct ways in which we can choose to
remove an upward-pointing triangle and downward-pointing triangle
from these six, according to whether the triangles are opposite or adjacent:
$$\cboard{
  AVAVAVA                 AVAVAVA
 AVAVAVAVA               AVAVAVAVA
AVAVA AVAVA             AVAV VAVAVA
VAVAV VAVAV             VAVA AVAVAV
 VAVAVAVAV               VAVAVAVAV
  VAVAVAV                 VAVAVAV
}
\vadjust{\vskip3pt}
$$
Such regions may be called ``holey hexagons'' of two different kinds.
Matt Blum tabulated the number of lozenge tilings of these regions,
for small values of $a=b=c$.
In the first (``opposite'') case, the number
of tilings of the holey hexagon is a nice round number (its greatest
prime factor appears to be bounded by a linear function of the size
of the region).  In the second (``adjacent'') case, the number
of tilings is not round.  Note, however, that in the second case,
the number of tilings of the holey hexagon divided by the number of tilings
of the unaltered hexagon (given to us by MacMahon's formula) is equal to
the probability that a random lozenge tiling of the hexagon contains a
lozenge that covers these two triangles; this probability tends to $\frac13$
for large $a$, at least on average \cite{PrCLP}.
Following this clue, we examine the difference between the aforementioned
probability (with its messy, un-round numerator) and the number $\frac13$.
The result is a fraction in which the numerator is now a nice round number.
So, in both cases, we have reason to think that there is an exact product
formula.

\begin{problem}
\label{prob4}
Determine the number of lozenge tilings of a
regular hexagon from which two of its innermost unit triangles
(one upward-pointing and one downward-pointing) have been removed.
\end{problem}

\begin{progress}
Theresia Eisenk\"olbl solved the first case of Problem 4
and Markus Fulmek and Christian Krattenthaler solved the second case.
Eisenk\"olbl \citeyear{PrE2} solves a generalization of the problem by applying
Mihai Ciucu's matchings factorization theorem, nonintersecting lattice
paths, and a nontrivial determinant evaluation.  Fulmek and Krattenthaler
\citeyear{PrFK2} compute the number of rhombus tilings of a hexagon
with sides $a,b,a,a,b,a$ (with $a$ and $b$ having the same parity)
that contain the rhombus that touches the center of the hexagon and
lies symmetric with respect to the symmetry axis that runs parallel
to the sides of length $b$.
For the proof of their formula they compute Hankel determinants
featuring Bernoulli numbers, which they do by using facts about continued
fractions, orthogonal polynomials, and, in particular, continuous Hahn
polynomials. The special case $a=b$ solves the second part of Problem 4.
\end{progress}

\hskip\parindent
I mentioned earlier that Kasteleyn's method, as interpreted by Percus,
allows one to write the number of matchings of a bipartite planar graph as
the determinant of a signed version of the bipartite adjacency matrix.
In the case of lozenge tilings of hexagons and the associated matchings,
it turns out that there is no need to modify signs of entries; the ordinary
bipartite adjacency matrix will do.  Greg Kuperberg \citeyear{PrKu2}
has noticed that when row-reduction and column-reduction are
systematically applied to the Kasteleyn--Percus matrix of
an $a,b,c$ semiregular hexagon,
one can obtain the $b$-by-$b$ Carlitz matrix \cite{PrCS}
whose $(i,j)$-th entry is
$a+c \choose a+i-j$.
(This matrix can also be recognized as the Gessel--Viennot matrix
that arises from interpreting each tiling
as a family of nonintersecting lattice paths \cite{PrGV}.)
Such reductions do not affect the determinant,
so we have a pleasing way of understanding the relationship
between the Kasteleyn--Percus matrix method
and the Gessel--Viennot lattice-path method.
In fact, such reductions
do not affect the \textit{cokernel} of the matrix
(an abelian group whose order is the determinant).
On the other hand, the cokernel of the Kasteleyn--Percus matrix
for the $a,b,c$ hexagon is clearly invariant under
permuting $a$, $b$, and $c$.
This gives rise to three different Carlitz matrices
that nontrivially have the same cokernel.
For example,
if $c=1$,
then one gets an $a$-by-$a$ matrix and a $b$-by-$b$ matrix
that both have the same cokernel,
whose structure can be determined ``by inspection'' if one notices
that the third Carlitz matrix of the trio is just a 1-by-1 matrix
whose sole entry is (plus or minus) a binomial coefficient.
In this special case,
the cokernel is just a cyclic group.

Greg Kuperberg poses this challenge:

\begin{problem}
\label{prob5}
Determine the cokernel of the Carlitz matrix,
or equivalently of the Kasteleyn--Percus matrix of the $a,b,c$ hexagon,
and if possible find a way to interpret the cokernel in terms of
the tilings.

This combines Questions 1 and 2 of Kuperberg \citeyear{PrKu2}.
As he points out in that article, in the case $a=b=c=2$,
one gets the noncyclic group $\Z/2\Z \times \Z/10\Z$ as the cokernel.
\end{problem}

\hskip\parindent
As was remarked above, one nice thing about the Kasteleyn--Percus matrices
of honeycomb graphs is that it is not necessary to make any of the
entries negative.  For general graphs,
however, there is no canonical way of defining $K$, in the sense
that there may be many ways of modifying the signs of certain
entries of the bipartite adjacency matrix of a graph so that all
nonzero contributions to the determinant have the same sign.
Thus, one should not expect the eigenvalues of $K$ to possess
combinatorial significance.  However, the spectrum of $K$ times
its adjoint $K^*$ is independent of which Kasteleyn--Percus matrix $K$
one chooses (as was independently shown by David Wilson and Horst
Sachs).  Thus, digressing somewhat from the topic of lozenge tilings,
we find it natural to ask:

\begin{problem}
\label{prob6}
What is the significance of the spectrum of
$K K^*$, where $K$ is any Kasteleyn--Percus matrix associated with a
bipartite planar graph?
\end{problem}

\begin{progress}
Nicolau Saldanha \citeyear{PrSal} has proposed a combinatorial
interpretation of the spectrum of $K K^*$.
Horst Sachs says (personal communication)
that $K K^*$ may have some significance
in the chemistry of polycyclic hydrocarbons (so-called
benzenoids) and related compounds as a useful approximate
measure of the ``degree of aromaticity''.
\end{progress}

\hskip\parindent
Returning now to lozenge tilings, or equivalently, matchings
of finite subgraphs of the infinite honeycomb,
consider the hexagon graph with $a=b=c=2$:%
$$
\vcenter{\rot{\psfig{file=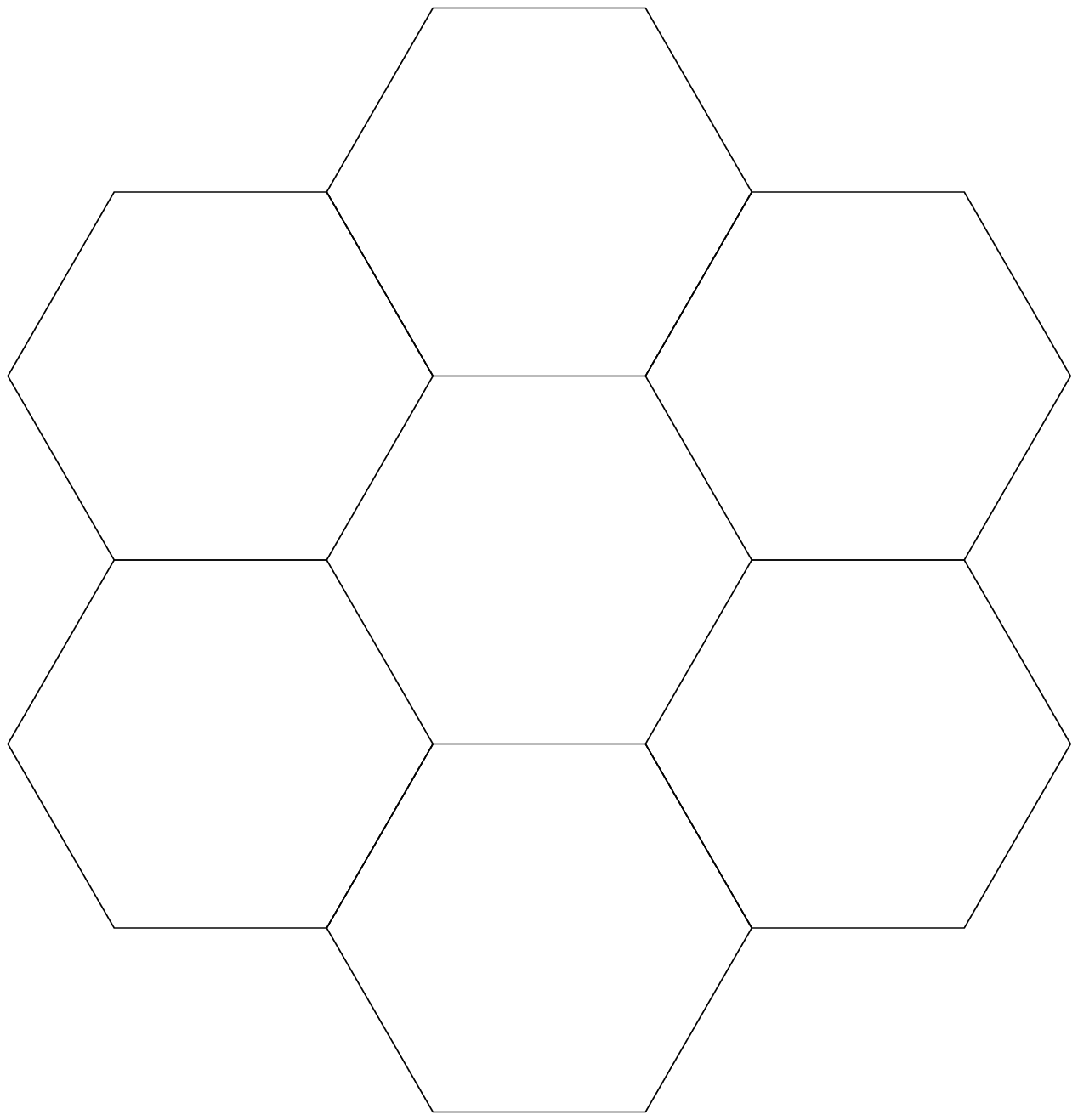,width=1.4in}}{0}{0}}
$$
This is the graph whose 20 matchings correspond to the 20 tilings
of the regular hexagon of side 2 by rhombuses of side 1.  If we look
at the probability of each individual vertical edge belonging to a
matching chosen uniformly at random (``edge-probabilities''), we get
\vadjust{\vskip4pt}
$$
\makeatletter
\vcenter{\rot{\psfig{file=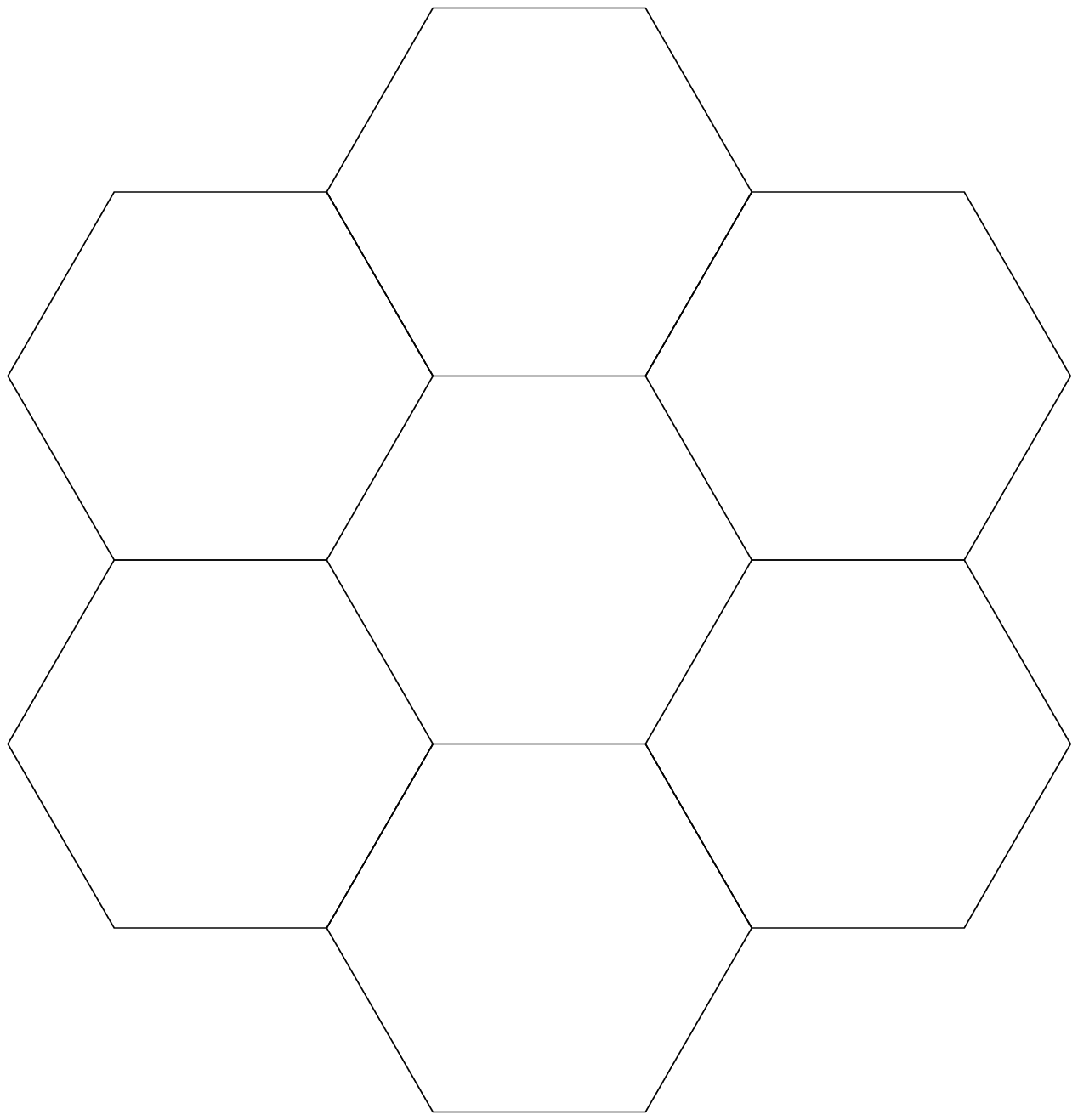,width=1.4in}}{0}{0}}
$$
Now look at this table of numbers as if it described a distribution
of mass.  If we assign the three rows $y$-coordinates $-1$ through 1,
we find that the weighted sum of the squares of the $y$-coordinates is equal to
$$(0.3+0.4+0.3)(-1)^2+(0.7+0.3+0.7+0.3)(0)^2+(0.3+0.4+0.3)(1)^2=2.$$
If we assign to the seven
columns $x$-coordinates $-3$ through 3, we find that the weighted sum of the
squares of the $x$-coordinates is equal to
$(0.7)(-3)^2+(0.6)(-2)^2+(0.3)(-1)^2+(0.8)(0)^2
+(0.3)(1)^2+(0.6)(2)^2+(0.7)(3)^2=20$.
\vadjust{\goodbreak}%
You can do a similar (but even easier) calculation yourself for the case
$a=b=c=1$, to see that the ``moments of inertia'' of the vertical
edge-probabilities around the horizontal and vertical axes are 0 and
1, respectively.  Using {\tt vaxmaple} to study the case $a=b=c=n$ for
larger values of $n$, I find that the moment of inertia about the horizontal
axis goes like
$$0, 2, 12, 40, 100, \dots$$
and the moment of inertia about the vertical axis goes like
$$1, 20, 93, 296, 725, \dots .$$
It is easy to show that the former moments of inertia
are given in general by the polynomial
$(n^4-n^2)/6$
(in fact, the number of vertical lozenges that have any particular
$y$-coordinate does not depend on the tiling chosen).
The latter moments of inertia are subtler;
they are not given by a polynomial of degree 4,
though it is noteworthy that the $n$-th term
is an integer divisible by $n$,
at least for the first few values of $n$.

\begin{problem}
\label{prob7}
Find the ``moments of inertia'' for the mass on
edges arising from edge-probabilities for random matchings of the
$a,b,c$ honeycomb graph.
\end{problem}

\section{Dominoes}

Now let us turn from lozenge-tiling problems to domino-tiling problems.
A \textit{domino} is a 1-by-2 or 2-by-1 rectangle.
Although lozenge tilings (in the guise of constrained plane partitions)
were studied first, it was really the study of domino tilings in Aztec
diamonds that gave current work on enumeration of matchings its current
impetus.  Here is the Aztec diamond of order 5:
$$\cboard{
    XX
   XXXX
  XXXXXX
 XXXXXXXX
XXXXXXXXXX
XXXXXXXXXX
 XXXXXXXX
  XXXXXX
   XXXX
    XX
}$$
A tiling of such a region by dominos is equivalent to a matching
of a certain (dual) subgraph of the infinite square graph.
This grid is bipartite, and it is convenient to color its vertices
alternately black and white;
equivalently, it is convenient to color the 1-by-1 squares alternately
black and white, so that every domino contains
one 1-by-1 square of each color.
Elkies, Kuperberg, Larsen, and Propp showed
in \cite{PrEKLP}
that the number of domino tilings
of such a region is $2^{n(n+1)/2}$ (where $2n$ is the number of rows),
and Gessel, Ionescu, and Propp proved in \cite{PrGIP} an exact formula
(originally conjectured by Jockusch) for the number of tilings of regions like
$$\cboard{
    XX
   XXXX
  XXXXXX
 XXXXXXXX
XXXX XXXXX
XXXX XXXXX
 XXXXXXXX
  XXXXXX
   XXXX
    XX
}$$
in which two innermost squares of opposite color have been removed.
(For some values of $n$, the number of tilings is exactly $\frac14$
times $2^{n(n+1)/2}$; in the other cases, there is an exact product
formula for the difference between the number of tilings and
$\bigl(\frac14\bigr)2^{n(n+1)/2}$.  It is this latter fact that motivated the
idea of trying something similar in the case of lozenge tilings,
as described in the paragraph preceding the statement of Problem 4.)

Now suppose one removes two squares from the middle of an Aztec diamond of
order $n$ in the following way:
$$\cboard{
    XX
   XXXX
  XXXXXX
 XXXX XXX
XXXXXXXXXX
XXXX XXXXX
 XXXXXXXX
  XXXXXX
   XXXX
    XX
}$$
(The two squares removed are a knight's-move apart, and subject to that
constraint, they are as close to being in the middle as they can be.
Up to symmetries of the square, there is only one way of doing this.)
The numbers of tilings one gets are as follows (for $n = 2$ through 10):
$$
\medmuskip 1mu
\displaylines{
2\cr
2^3\cr
2^5\cdot 5\cr
2^9\cdot 3^2\cr
2^{17}\cdot 3\cr
2^{22}\cdot 3^2\cr
2^{24}\cdot 3^2\cdot73\cr
2^{31}\cdot 3^2\cdot5^2\cdot11\cr
2^{47}\cdot 3^2\cdot5\cr
}
$$
Only the presence of the large prime factor 73 makes one doubt
that there is a general product formula; the other prime factors are
reassuringly small.

\begin{problem}
\label{prob8}
Count the domino tilings of an Aztec diamond from which
two close-to-central squares, related by a knight's move, have been
deleted.
\end{problem}

\begin{progress}
Harald Helfgott has solved this problem;
it follows from the main result in his thesis \citeyear{PrH}.
The formula is somewhat complicated,
as the prime factor 73 might have led us to expect.
(One of the factors in Helfgott's product formula
is a single-indexed sum;
73 arises as $128-60+5$.)
\end{progress}

\hskip\parindent
One can also look at ``Aztec rectangles'' from which squares have been
removed so as to restore the balance between black and white squares
(a necessary condition for tileability).  For instance, one can remove
the central square from an $a$-by-$b$ Aztec rectangle in which $a$
and $b$ differ by 1, with the larger of $a,b$ odd:
$$\cboard{
   XX
  XXXX
 XXXXXX
XXXXXXXX
XXXX XXXX
 XXXXXXXX
  XXXXXX
   XXXX
    XX
}$$

\begin{problem}
\label{prob9}
Find a formula for the number of domino tilings of a
$2n$-by-$(2n+1)$ Aztec rectangle with its central square removed.
\end{problem}

\begin{progress}
This had already been solved when I posed the problem; it is
a special case of a result of Ciucu \citeyear[Theorem 4.1]{PrCi1}.
Eric Kuo solved the problem independently.
\end{progress}

\hskip\parindent
What about $(2n-1)$-by-$2n$ rectangles?  For these regions, removing
the central square does not make the region tilable.  However, if
one removes any one of the four squares adjacent to the middle square,
one obtains a region that is tilable, and moreover, for this region
the number of tilings appears to be a nice round number.

\begin{problem}
\label{prob10}
Find a formula for the number of domino tilings of a
$(2n-1)$-by-$2n$ Aztec rectangle with a square adjoining the central
square removed.
\end{problem}

\begin{progress}
This problem was solved independently three times:\,\
by Harald Helfgott and Ira Gessel \citeyear{PrHG},
by Christian Krattenthaler \citeyear{PrKr},
and by Eric Kuo (private communication).
Gessel and Helfgott solve a more general problem than Problem 10.
Krattenthaler's preprint
gives several results concerning the enumeration of matchings
of Aztec rectangles where (a suitable number of) collinear vertices are
removed, of which Problem 10 is just a special case.
There is some overlap between the results of Helfgott and Gessel
and the results of Krattenthaler.
\end{progress}

\hskip\parindent
At this point, some readers may be
wondering why $m$-by-$n$ rectangles have not played
a bigger part in the story.
Indeed, one of the surprising facts of life in the study of enumeration
of matchings is that Aztec diamonds and their kin have been
much more fertile ground for exact combinatorics that the seemingly
more natural rectangles.  There are, however,
a few cases I know of in which something rather nice turns up.  One
is the problem of Ira Gessel that appears as Problem 20 in this document.
Another is the work done by Jockusch \citeyear{PrJ} and, later, Ciucu \citeyear{PrCi1}
on why the number
of domino tilings of the square is always either a perfect square or
twice a perfect square.  In the spirit of the work of Jockusch and
Ciucu, I offer here a problem based on Lior Pachter's observation
\cite{PrPK} that the region on the left below, obtained by removing
8 dominos from a 16-by-16 square, has exactly one tiling.
What if we make the intrusion half as long, as in the region on the right?
$$\cboard{
XXXXXXXXXXXXXXXX
XXXXXXXXXXXXXXXX
XXXXXXXXXXXXXXXX
XXXXXXXXXXXXXXXX
XXXXXXXXXXXXXXXX
XXXXXXXXXXXXXXXX
XXXXXXXXXXXXXXXX
XXXXXXXXXXXXXXXX
XXXXXXX  XXXXXXX
XXXXXX  XXXXXXXX
XXXXX  XXXXXXXXX
XXXX  XXXXXXXXXX
XXX  XXXXXXXXXXX
XX  XXXXXXXXXXXX
X  XXXXXXXXXXXXX
  XXXXXXXXXXXXXX
}\hskip .25\hsize
\cboard{
XXXXXXXXXXXXXXXX
XXXXXXXXXXXXXXXX
XXXXXXXXXXXXXXXX
XXXXXXXXXXXXXXXX
XXXXXXXXXXXXXXXX
XXXXXXXXXXXXXXXX
XXXXXXXXXXXXXXXX
XXXXXXXXXXXXXXXX
XXXXXXXXXXXXXXXX
XXXXXXXXXXXXXXXX
XXXXXXXXXXXXXXXX
XXXXXXXXXXXXXXXX
XXX  XXXXXXXXXXX
XX  XXXXXXXXXXXX
X  XXXXXXXXXXXXX
  XXXXXXXXXXXXXX
}$$
That is, we take a $2n$-by-$2n$ square (with $n$ even) and remove $n/2$
dominos from it, in a partial zig-zag pattern that starts from the corner.
Here are the numbers we get, in factored form, for $n=2,4,6,8,10$:
$$
\medmuskip 1mu
\displaylines{
2\cdot3^2\cr
2^2\cdot3^6\cdot13^2\cr
2^3\cdot3^2\cdot5^4\cdot7^2\cdot3187^2\cr
2^4\cdot11771899^2\cdot27487^2\cr
2^5\cdot2534588575976069659^2
}$$
The factors are ugly, but the exponents are nice: we get $2^{n/2}$ times
an odd square.

Perhaps this is a special case of a two-parameter fact that says that
you can take an intrusion of length $m$ in a $2n$-by-$2n$ square and the
number of tilings of the resulting region will always be a square or
twice a square.

\begin{problem}
\label{prob11}
What is going on with ``intruded Aztec diamonds''?
In particular, why is the number of tilings so square-ish?
\end{problem}

\hskip\parindent
It should also be noted that the square root of the odd parts of these
numbers (3, $3^3 \cdot 13$, etc.)\ alternate between 1 and 3 mod 4.
Perhaps these quantities are continuous functions of $n$ in the 2-adic
sense, as is the case for intact $2n$-by-$2n$ squares \cite{PrCo};
however, the presence of large prime factors means that no simple
product formula is available, and that the analysis will require new
techniques.

\medbreak

We now return to the Kasteleyn--Percus matrices discussed earlier.
Work of Rick Kenyon and David Wilson  \cite{PrKe}
has shown that the \textit{inverses}
of these matrices are loaded with combinatorial information, so it
would be nice to get our hands on them.  Unfortunately, there are
many nonzero entries in the inverse-matrices.  (Recall that
the Kasteleyn--Percus matrices themselves, being nothing more than adjacency
matrices in which some of the 1's have been strategically replaced
by $-1$'s, are sparse; their inverses, however, tend to have most
if not all of their entries nonzero.)  Nonetheless, some exploratory
``numerology'' leaves room for hope that this is do-able.

Consider the Kasteleyn--Percus matrix $K_n$ for the Aztec diamond of
order $n$, in which every vertical domino with its white square on top
(relative to some fixed checkerboard coloring) has its sign
inverted\emdash that is, the corresponding 1 in the bipartite
adjacency matrix is replaced by $-1$.

\begin{problem}
\label{prob12}
Show that the sum of the entries of the matrix inverse
of $K_n$ is $\frac12(n-1)(n+3) - 2^{n-1} + 2$.
\end{problem}

\hskip\parindent
(This formula works for $n=1$ through $n=8$.)

\begin{progress}
Harald Helfgott has solved a similar problem using the main result of
his thesis
\citeyear{PrH}, and it is likely that the result asserted in Problem 12
can be proved similarly.
(A slight technical hurdle arises from the fact that Helfgott's thesis uses
a different sign-convention for the Kasteleyn--Percus matrix, which
results in different signs, and a different sum, for the inverse
matrix; however, Helfgott's methods are quite general, so there is
no conceptual obstacle to applying them to Problem 12.)

I should mention that my original reason for examining the
sum of the entries of the inverse Kasteleyn--Percus matrix
was to see whether there might be formulas governing the
individual entries themselves.  Helfgott's work provides
such formulas.

Also, in this connection, Greg Kuperberg and Douglas
Zare have some high-tech ruminations on the inverses of Kasteleyn--Percus
matrices, and there is a chance that representation-theory methods will
give a different way of proving the result.
\end{progress}

\hskip\parindent
Now we turn to a class of regions I call ``pillows''.
Here are a ``0 mod 4'' pillow of ``order 5'' and a
``2 mod 4'' pillow of ``order 7'':
$$\cboard{
            XXXX
         XXXXXXXX
      XXXXXXXXXXXX
   XXXXXXXXXXXXXXXX
XXXXXXXXXXXXXXXXXXXX
XXXXXXXXXXXXXXXXXXXX
 XXXXXXXXXXXXXXXX
  XXXXXXXXXXXX
   XXXXXXXX
    XXXX
}\hskip.1\hsize
\cboard{
                  XX
               XXXXXX
            XXXXXXXXXX
         XXXXXXXXXXXXXX
      XXXXXXXXXXXXXXXXXX
   XXXXXXXXXXXXXXXXXXXXXX
XXXXXXXXXXXXXXXXXXXXXXXXXX
XXXXXXXXXXXXXXXXXXXXXXXXXX
 XXXXXXXXXXXXXXXXXXXXXX
  XXXXXXXXXXXXXXXXXX
   XXXXXXXXXXXXXX
    XXXXXXXXXX
     XXXXXX
      XX
}$$
It turns out (empirically) that the number of tilings of the 0-mod-4 pillow
of order $n$ is a perfect square times the coefficient of $x^n$ in the Taylor
expansion of $(5+3x+x^2-x^3)/(1-2x-2x^2-2x^3+x^4)$.  This fact came to light
in several steps.  First it was noticed that the number of tilings has a
comparatively small square-free part.  Then it was noticed that in the
derived sequence of square-free parts, many terms were roughly three times
the preceding term.  Then it was noticed that, by 
judiciously including some of the square factors, one could obtain
a sequence in which each term was roughly three times the preceding term.
Finally it was noticed that this approximately geometric sequence satisfied
a fourth-order linear recurrence relation.

Similarly, it appears that the number of tilings of the 2-mod-4 pillow of
order $n$ is a perfect square times the coefficient of $x^n$ in the Taylor
expansion of $(5+6x+3x^2-2x^3)/(1-2x-2x^2-2x^3+x^4)$.  (If you are wondering
about ``odd pillows'', I should mention that there is a nice formula for the
number of tilings, but this is not an interesting result, because an odd pillow
splits up into many small noncommunicating sub-regions such that a tiling
of the whole region corresponds to a choice of tiling on each of the
sub-regions.)

\begin{problem}
\label{prob13}
Find a general formula for the number of domino tilings
of even pillows.
\end{problem}

\hskip\parindent
Jockusch looked at the Aztec diamond of order $n$ with a 2-by-2
hole in the center, for small values of $n$;
he came up with a conjecture for the number of domino tilings,
subsequently proved by Gessel, Ionescu, and Propp \cite{PrGIP}.
One way to generalize this is to make the hole larger,
as was suggested by Douglas Zare
and investigated by David Wilson.
Here is an abridged and adapted version of the report
David Wilson sent me on October 15, 1996:

\medskip

Define the Aztec window with outer order $y$ and
inner order $x$ to be the Aztec diamond of order $y$ with an
Aztec diamond of order $x$ deleted from its center.  For
example, this is the Aztec window with orders 8 and 2:
$$\cboard{
       XX
      XXXX
     XXXXXX
    XXXXXXXX
   XXXXXXXXXX
  XXXXXXXXXXXX
 XXXXXX  XXXXXX
XXXXXX    XXXXXX
XXXXXX    XXXXXX
 XXXXXX  XXXXXX
  XXXXXXXXXXXX
   XXXXXXXXXX
    XXXXXXXX
     XXXXXX
      XXXX
       XX
}$$

There are a number of interesting patterns that show up
when we count tilings of Aztec windows.  For one thing,
if $w$ is a fixed even number, and $y = x+w$, then for any $w$
the number of tilings appears to be a polynomial in $x$.
(When $w$ is odd, and $x$ is large enough, there are no tilings.)
For $w=6$, the polynomial is
$$
\displaylines{
 8192 x^8  + 98304 x^7  + 573440 x^6  + 2064384 x^5  + 4988928 x^4 \hfill\cr
\hfill + 8257536 x^3  + 9175040 x^2  + 6291456 x + 2097152. }
$$
This can be written as
$$
2^{17} \bigl( \tfrac{1}{2} \bigl(x+\tfrac{3}{2}\bigr)^2 + \tfrac{7}{8} \bigr)^4
$$
or as
$$
2^{17}x^4 \ \ \circ \ \
\tfrac{1}{2}x+\tfrac{7}{8} \ \ \circ \ \ \bigl(x+\tfrac{3}{2}\bigr)^2,
$$
where it is understood that these three polynomials get composed.

More generally, all the polynomials in $x$ that arise in this fashion
appear to ``factor'' in the sense of functional composition.  Here
are the factored forms of the polynomials for $n=2,4,6,8,10$:
$$
\arraycolsep=4pt
\begin{array}{rcccl}
2^{3}x^4&\circ&1&\circ&\bigl(x+\frac12\bigr)^2\\[2ex]
2^{8}x^2&\circ&x+1&\circ&\bigl(x+1\bigr)^2\\[2ex]
2^{17}x^4&\circ&\frac12x+\frac{7}{8}&\circ&\bigl(x+\frac32\bigr)^2\\[2ex]
2^{28}x^2&\circ&\frac{1}{144}x^4+\frac{7}{72}x^3
+\frac{41}{144}x^2+\frac{11}{18}x+1&\circ&\bigl(x+2\bigr)^2\\[2ex]
2^{43}x^4&\circ&\frac{1}{144}x^3+\frac{61}{576}x^2
+\frac{451}{2304}x+\frac{967}{1024}&\circ&\bigl(x+\frac52\bigr)^2
\end{array}
$$
In general the rightmost polynomial is $(x+w/4)^2$,
and the leftmost polynomial is either a perfect square,
twice a fourth power, or half a fourth power, depending
on $w$ mod 8.  A pattern for the middle polynomial however
is elusive.

\begin{problem}
\label{prob14}
Find a general formula for the number of domino tilings
of Aztec windows.
\end{problem}

\begin{progress}
Constantin Chiscanu
found a polynomial bound on
the number of domino tilings of the Aztec window
of inner order $x$ and outer order $x+w$ \cite{PrCh}.
Douglas Zare used the transfer-matrix method to show that
the number of tilings is not just bounded by a polynomial,
but given by a polynomial, for each fixed $w$ \cite{PrZ}.
\end{progress}

\section{Miscellaneous}

Now we come to some problems involving tiling that fit neither the
domino-tiling nor the lozenge-tiling framework.
Here the more general picture is that we have some periodic dissection
of the plane by polygons,
such that an even number of polygons meet at each vertex,
allowing us to color the polygons alternately black or white.
We then make a suitable choice of a finite region $R$
composed of equal numbers of black and white polygons,
and we look at the number of ``diform'' tilings of the region,
where a \textit{diform} is the union of two polygonal cells
that share an edge.
In the case of domino tilings, the underlying dissection of the infinite plane
is the tiling by squares, 4 around each vertex;
in the case of lozenge tilings, the underlying dissection of the infinite plane
is the tiling by equilateral triangles, 6 around each vertex.

Other sorts of periodic dissections have already played a role in
the theory of enumeration of matchings.
For instance, there is a tiling of the plane by isosceles right triangles
associated with a discrete reflection group in the plane;
in this case, the right choice of $R$ (see Figure 3)
gives us a region that can be tiled in $5^{n^2/4}$ ways when $n$ is even
and in $5^{(n^2-1)/2}$ or $2 \cdot 5^{(n^2-1)/2}$ ways when $n$ is odd
\cite{PrY}.
\begin{figure}
\vskip1pt
\centerline{\psfig{file=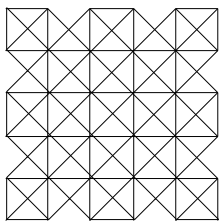}}
\vskip-1pt
\caption{A fortress of order 5, with $2 \times 5^6$ diform tilings.}
\label{fortress}
\end{figure}

\begin{figure}[b]
\vskip1pt
\centerline{\psfig{file=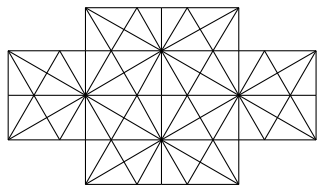}}
\vskip-1pt
\caption{An Aztec dungeon of order 2, with $13^3$ diform tilings.}
\label{azdungeon}
\end{figure}

Similarly, in the tiling of the plane by triangles
that comes from a 30 degree, 60 degree, 90 degree right triangle
by repeatedly reflecting it in its edges,
a certain region called the ``Aztec dungeon''
(see Figure 4)
gives rise to a tiling problem in which powers of 13 occur
(as was proved in not-yet-published work of Mihai Ciucu).

A key feature of these regions $R$
is revealed by looking at the colors of those polygons in the dissection
that share an edge with the border of $R$.
One sees that the border splits up into four long stretches
such that along each stretch,
all the polygons that touch the border have the same color.
It is not clear why regions with this sort of property
should be the ones that give rise to
the nicest enumerations,
but this appears to happen in practice.

One interesting case
arises from a rather symmetric
dissection of the plane into
equilateral triangles, squares, and regular hexagons,
with 4 polygons meeting at each vertex
and with no two squares sharing an edge.
A typical diform tiling of this region
(called a ``dragon'')
is shown in Figure 5.
\begin{figure}
\centerline{\psfig{file=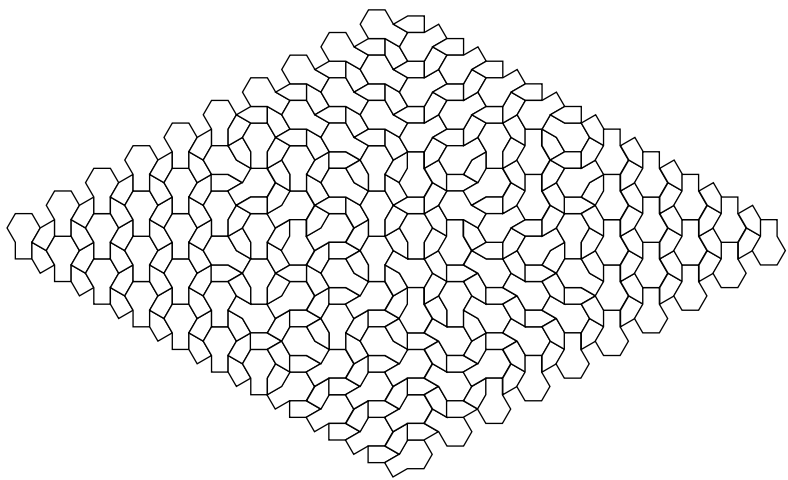}}
\caption{A dragon of order 10 (tiled).}
\label{trisquex}
\end{figure}
Empirically, one finds that the number of diform tilings is $2^{n(n+1)}$.

\begin{problem}
\label{prob15}
Prove that the number of diform
tilings of the dragon of order $n$ is $2^{n(n+1)}$.
\end{problem}

\begin{progress}
Ben Wieland solved this problem (private communication).
\end{progress}

\hskip\parindent
Incidentally, the tiling shown in Figure 5 was generated
using an algorithm that generates each of the possible diform tilings
of the region with equal probability.
It is no fluke that the tiling looks so orderly
in the left and right corners of the region;
this appears to be typical behavior in situations of this kind.
This phenomenon has been analyzed rigorously for two tiling-models:
lozenge tilings of hexagons \cite{PrCLP} and domino tilings of Aztec diamonds
\cite{PrCEP}.

One way to get a new dissection of the plane from an old one is to refine it.
For instance, starting from the dissection of the plane into squares,
one can draw in every $k$-th southwest-to-northeast diagonal.  When
$k$ is 1, this is just a distortion of the dissection of the plane into
equilateral triangles.  When $k$ is 2, this is a dissection that leads to
finite regions for which the number of diform tilings is a known power of 2,
thanks to a theorem of Chris Douglas \citeyear{PrD}.
But what about $k=3$ and higher?

For instance, we have the roughly hexagonal region shown in Figure 6;
certain boundary vertices have been marked with a dot so as to
bring out the large-scale $2,3,2,2,3,2$ hexagonal structure more clearly.

\begin{figure}[b]
\centerline{\psfig{file=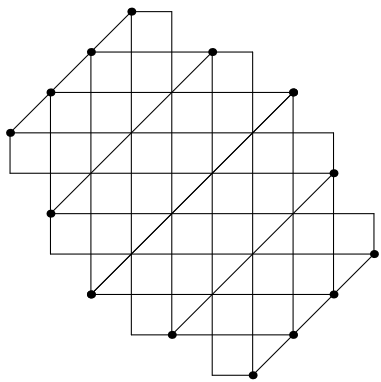,width=2in}}
\caption{A region for Problem 16.}
\label{wieland}
\end{figure}

The cells of this region are triangles and squares.
The region has $17920 = 2^9 \cdot 5 \cdot 7$ diform tilings.

\begin{problem}
\label{prob16}
Find a formula for the number of diform tilings in the
$a,b,c$ quasihexagon in the dissection of the plane that arises from
slicing the dissection into squares along every third upward-sloping
diagonal.
\end{problem}

\hskip\parindent
One reason for my special interest in Problem 16
is that it seems to be a genuine hybrid of domino tilings of Aztec diamonds
and lozenge tilings of hexagons.

\begin{progress}
Ben Wieland solved this problem in the case $a=b=c$ (which, as it
turns out, is also the solution to the case $a=b<c$ and the case
$a=c<b$).  In these cases the number of tilings is always a power
of two.  The general case does not yield round numbers, so there
is no simple product formula.
\end{progress}

\hskip\parindent
The approach underlying Ben Wieland's solutions to the last two problems
is a method of subgraph substitution that has already been of
great use in enumeration of matchings of graphs.
I will not go into great detail here on this method
[\citeNP{PrPr1}; \citeyearNP{PrPr2}],
but here is an overview:
One studies graphs with weights assigned to their edges,
and one does weighted enumeration of matchings,
where the weight of a matching is the product
of the weights of the constituent edges.
One then looks at local substitutions of subgraphs within a graph that
preserve the sum of the weights of the matchings,
or more generally,
multiply the sum of the weights of the matchings
by some predictable factor.
Then the problem of weight-enumerating matchings of one graph
reduces to the problem of weight-enumerating matchings
of another graph.
Iterating this procedure,
one can often eventually reduce the graph
to something easier to understand.

\medskip

Problems 15 and 16 are just two instances of a broad class of problems
arising from periodic graphs in the plane.  A unified understanding
of this class of problems has begun to emerge, by way of subgraph
substitution.  The most important open problem connected with this class
of results is the following:

\begin{problem}
\label{prob17}
Characterize those local substitutions that have
a predictable effect on the weighted sum of matchings of a graph.
\end{problem}

\begin{figure}[b]
\centerline{\psfig{file=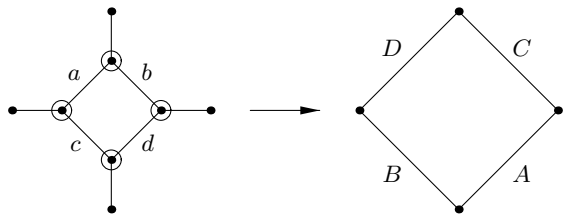,width=3in}}
\caption{The ``urban renewal'' substitution.}
\label{urban}
\end{figure}

\hskip\parindent 
The most useful local substitution so far has been the one shown in
Figure~7, where unmarked edges have weight 1 and where $A,B,C,D$ are
respectively obtained from $a,b,c,d$ by dividing by $ad+bc$; if $G$
and $G'$ denote the graph before and after the substitution, one can
check that the sum of the weights of the matchings of $G'$ equals the
sum of the weights of the matchings of $G$ divided by $ad+bc$.

It is required that the four innermost vertices
have no neighbors other than the four vertices shown;
this constraint is indicated by circling them.
Noncircled vertices may have any number of neighbors.

\begin{figure}[h]
\vskip2pt
\centerline{\psfig{file=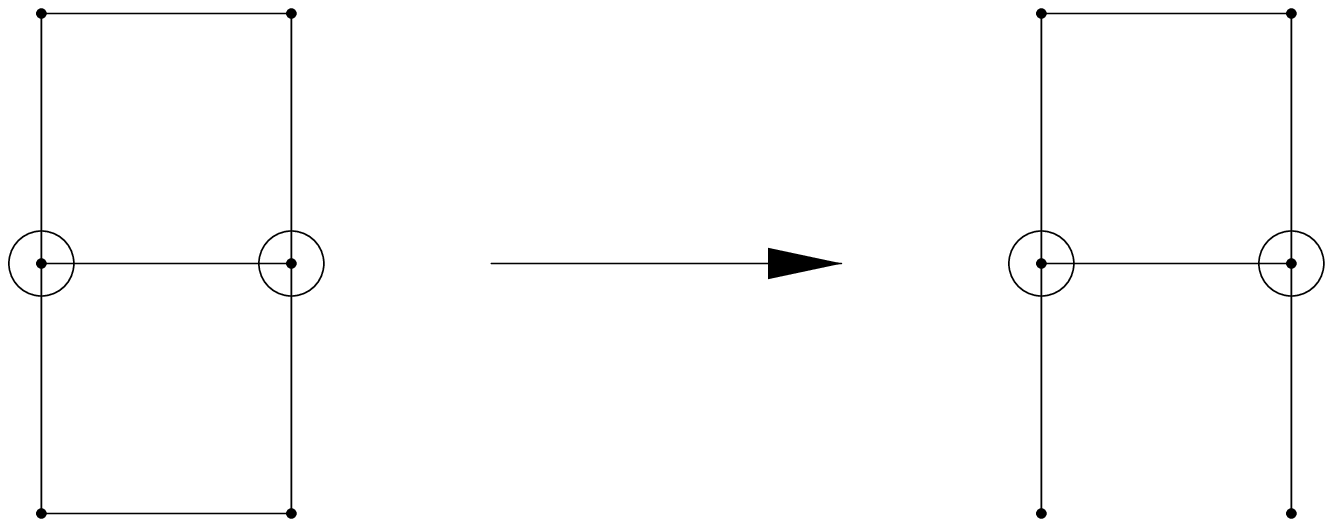,width=2in}}
\vskip-2pt
\caption{Rick Kenyon's substitution.}
\label{kenyon}
\end{figure}

The substitution shown in Figure 8
(a straightforward generalization of
a clever substitution due to Rick Kenyon)
has also been of use.
Here the new weights are not entirely determined by the old,
but have a single degree of freedom;
the relevant formulas can be written as
$$
\def\quad{\hskip 10 pt minus 6pt}
A  =  \frac{abc+aeg+cdf}{bc+eg} \,, \quad B  =  b \,, \quad
D  =  \frac{dg}{bc+eg} E \,, \quad
F  =  ef \frac{1}{E} \,,\quad
G  =  (bc+eg) \frac{1}{E}\,,
$$
with $E$ free.
As before, the circled vertices
must not have any neighbors
other than the ones shown.
In this case, the sum of the weights in the before-graph $G$
is exactly equal to the sum of the weights in the after-graph $G'$;
there is no need for a correction factor
like the $1/(ad+bc)$ that arises in urban renewal.

The extremely powerful ``wye-delta'' substitution of
Colbourn, Provan, and Vertigan \cite{PrCPV} should also be mentioned.

\medskip

Up till now we have been dealing exclusively with bipartite planar graphs.
We now turn to the less well-explored nonbipartite case.

For instance, one can look at the triangle graph of order $n$,
shown in Figure~9 in the case $n=4$.
(Here $n$ is the number of vertices in the longest row.)

\begin{figure}
\centerline{\psfig{file=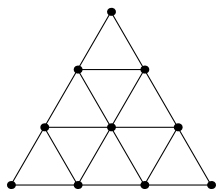,width=1.1in}}
\caption{The triangle graph.}
\label{triangle}
\end{figure}

Let $M(n)$ denote the number of matchings of the triangle graph of order $n$.
When $n$ is 1 or 2 mod 4, the graph has an odd number of vertices
and $M(n)$ is 0; hence let us only consider the cases in which
$n$ is 0 or 3 mod 4.
Here are the
first few values of $M(n)$, expressed in factored form:
$2$,
$2 \cdot 3$,
$2 \cdot 2 \cdot 3 \cdot 3 \cdot 61$,
$2 \cdot 2 \cdot 11 \cdot 29 \cdot 29$,
$2^3 \cdot 3^3 \cdot 5^2 \cdot 7^2 \cdot 19 \cdot 461$,
$2^3 \cdot 5^2 \cdot 37^2 \cdot 41 \cdot 139^2$,
$2^4 \cdot 73 \cdot 149 \cdot 757 \cdot 33721 \cdot 523657$,
$2^4 \cdot 3^8 \cdot 17 \cdot 37^2 \cdot 703459^2$,
\dots.
It is interesting that $M(n)$ seems to be divisible by
$2^{\lfloor (n+1)/4 \rfloor}$
but no higher power of 2;
it is also interesting that when we divide by this power of 2,
in the case where $n$ is a multiple of 4,
the quotient we get, in addition to being odd,
is a perfect square times a small number
$(3, 11, 41, 17, \dots{})$.

\begin{problem}
\label{prob18}
How many matchings does the triangle graph
of order $n$ have?
\end{problem}

\begin{progress}
Horst Sachs \citeyear{PrSac2} has responded to this problem.
\end{progress}

\hskip\parindent
One can also look at graphs that are bipartite but not planar.
A natural example is the $n$-cube (that is, the $n$-dimensional cube
with $2^n$ vertices).  It has been shown that the number
of matchings of the $n$-cube goes like $1$, $2$, $9=3^2$,
$272=16 \cdot 17$,
$589185=3^2 \cdot 5 \cdot 13093$, \dots.

\begin{problem}
\label{prob19}
Find a formula for the number of matchings
of the $n$-cube.
\end{problem}

\hskip\parindent
(This may be intractable; after all, the graph has exponentially
many vertices.)

\begin{progress}
L\'aszl\'o Lov\'asz gave a simple proof of my (oral) conjecture that
the number of matchings of the $n$-cube
has the same parity as $n$ itself.
Consider the orbit of a particular matching of the $n$-cube
under the group generated by the $n$ standard reflections of the $n$-cube.
If all the edges are parallel (which can happen in exactly $n$ ways),
the orbit has size 1;
otherwise the size of the orbit is
of the form $2^k$ (with $k \geq 1$)\emdash an even number.
The claim follows,
and similar albeit more complex reasoning
should allow one to compute the enumerating sequence
modulo any power of 2.
Meanwhile,
L. H.\ Clark, J. C.\ George, and T. D.\ Porter have shown \cite{PrCGP}
that if one lets $f(n)$ denote the number of 1-factors in the $n$-cube,
then $$f(n)^{2^{1-n}} \sim n/e$$
as $n \rightarrow \infty$.
It was subsequently pointed out by Bruce Sagan
that the main result of Clark et al.\ \citeyear{PrCGP} is a special case
of the theorem cited by Lov\'asz and Plummer \citeyear[top of page 312]{PrLP}.
\end{progress}

\hskip\parindent
Finally, we turn to a problem involving domino tilings of rectangles,
submitted by Ira Gessel (what follows are his words):

%[Note to people reading this over the Web: if you are using {\tt ghostview},
%you will have to switch over to {\tt gs} to preview pages 24 and 25 of this
%document.  If you're good at Postscript and can think of a way for me to
%fix the problem, please let me know!]

%Make ! a letter so we can use it in the names of control sequences
\catcode `!=11

\newdimen\squaresize
\newdimen\thickness
\newdimen\Thickness
\newdimen\ll! \newdimen \uu! \newdimen\dd!
\newdimen \rr! \newdimen \temp!

%parameters are left, up, down, right, and contents
\def\sq!#1#2#3#4#5{%
\ll!=#1 \uu!=#2 \dd!=#3 \rr!=#4
\setbox0=\hbox{%
%left edge
 \temp!=\squaresize\advance\temp! by .5\uu!
 \rlap{\kern -.5\ll!
 \vbox{\hrule height \temp! width#1 depth .5\dd!}}%
%
%up edge
 \temp!=\squaresize\advance\temp! by -.5\uu!
 \rlap{\raise\temp!
 \vbox{\hrule height #2 width \squaresize}}%
%
%down edge
 \rlap{\raise -.5\dd!
 \vbox{\hrule height #3 width \squaresize}}%
%
%right edge
 \temp!=\squaresize\advance\temp! by .5\uu!
 \rlap{\kern \squaresize \kern-.5\rr!
 \vbox{\hrule height \temp! width#4 depth .5\dd!}}%
%
%contents
 \rlap{\kern .5\squaresize\raise .5\squaresize
 \vbox to 0pt{\vss\hbox to 0pt{\hss $#5$\hss}\vss}}%
}%end of \hbox
 \ht0=0pt \dp0=0pt \box0
}%end of \sq!

\def\vsq!#1#2#3#4#5\endvsq!{\vbox
  to \squaresize{\hrule width\squaresize height 0pt%
\vss\sq!{#1}{#2}{#3}{#4}{#5}}}

\newdimen \LL! \newdimen \UU! \newdimen \DD! \newdimen \RR!

\def\vvsq!{\futurelet\next\vvvsq!}
\def\vvvsq!{\relax
  \ifx     \next l\LL!=\Thickness \let\continue!=\skipnexttoken!
  \else\ifx\next u\UU!=\Thickness \let\continue!=\skipnexttoken!
  \else\ifx\next d\DD!=\Thickness \let\continue!=\skipnexttoken!
  \else\ifx\next r\RR!=\Thickness \let\continue!=\skipnexttoken!
  \else\ifx\next P\let\continue!=\place!
  \else\def\continue!{\vsq!\LL!\UU!\DD!\RR!}%
  \fi\fi\fi\fi\fi
  \continue!}

\def\skipnexttoken!#1{\vvsq!}

\def\place! P#1#2#3{%
\rlap{\kern.5\squaresize\temp!=.5\squaresize\kern#1\temp!
  \temp!=\squaresize
  \advance\temp! by #2\squaresize \temp!=.5\temp!
  \raise\temp!\vbox
   to 0pt{\vss\hbox to 0pt{\hss$#3$\hss}\vss}}\vvsq!}

\def\Young#1{\LL!=\thickness \UU!=\thickness
 \DD! = \thickness \RR! = \thickness
\vbox{\smallskip\offinterlineskip
\halign{&\vvsq! ## \endvsq!\cr #1}}}

\def\Youngt#1{\LL!=\thickness \UU!=
  \thickness \DD! = \thickness \RR! = \thickness
\vtop{\offinterlineskip
\halign{&\vvsq! ## \endvsq!\cr #1}}}

\def\blank{\omit\hskip\squaresize}
\catcode `!=12

%\magnification1200
\thickness=.4pt

We consider dimer coverings of an $m\times n$ rectangle,
with $m$ and $n$ even. We assign a vertical domino from row
$i$ to row $i+1$ the weight
$\sqrt {y_i}$ and a horizontal domino from column $j$ to
column $j+1$ the weight $\sqrt {x_j}$. For example, the
covering
$$\Thickness=0pt
\squaresize=30pt
\Young{dP0{-1}{\sqrt{y_1}}&rP10
  {\sqrt{x_2}}&l&dP0{-1}
  {\sqrt{y_1}}&
  rP10{\sqrt{x_5}}&l&rP10{\sqrt{x_7}}
  &l&dP0{-1}{\sqrt{y_1}}&dP0{-1}{\sqrt{y_1}}\cr
  u&rP10{\sqrt{x_2}}&l&u&rP10
  {\sqrt{x_5}}&l&rP10{\sqrt{x_7}}&l&u&u\cr}
$$
for $m=2$ and $n=10$ has weight $y_1^2
x_2x_5x_7$. (The weight will always be a product
of  integral powers of the $x_i$ and $y_j$.)

Now I'll define what I call ``dimer tableaux.''
Take an $m/2$ by $n/2$ rectangle and split it
into two parts by a path from the lower left
corner to the upper right corner. For example
(with $m=6$ and $n=10$)
$$\Thickness=1pt
\squaresize 15pt
\Young{&&&&r\cr
&&&&lu\cr
d&dr&u&u&\cr}
$$
Then fill in the upper left
part with entries from 1, 2, \dots, $n-1$ so that
for adjacent entries
$\squaresize 10pt\lower 2pt\vbox{\yij}$
we have $i<j-1$ and for adjacent entries
$\squaresize 10pt\lower 8pt\vbox{\yIJ}$ we have $i\le j+1$, and fill in the
lower-right partition with entries from
$1,2\ldots, m-1$ with the reverse inequalities (
$\squaresize 10pt\lower 2pt\vbox{\yij}$
implies  $i\le j+1$ and
$\squaresize 10pt\lower 8pt\vbox{\yIJ}$
implies  $i<j-1$). We weight an $i$ in
the upper-left part by $x_i$ and a $j$ in the
lower-right part by $y_j$.

\begin{theorem}
The sum of the
weights of the $m\times n$ dimer coverings is
equal to the sum of the weights of the $m/2\times
n/2$ dimer tableaux.
\end{theorem}

My proof is not very enlightening; it essentially
involves showing that both of these are counted
by the same formula.

\begin{problem}
\label{prob20}
Is there an ``explanation'' for
this equality? In particular, is there  a
reasonable bijective proof?  Notes:
\begin{enumerate}
\item[(1)] The case $m=2$ is easy: the $2\times
10$ dimer covering above corresponds to the
$1\times 5$ dimer tableau
$$\squaresize
20pt\Thickness1pt\Young{dx_2&dx_5&drx_7&uy_1&uy_1\cr}$$
(there's only one possibility!).
\item[(2)] If we set $x_i=y_i=0$ when $i$ is even
(so that every two-by-two square of the dimer
covering may be chosen independently), then the
equality is equivalent to the identity
$$\prod_{i,j}(x_i+y_j)=\sum_{\lambda}s_{\lambda}(x)s_{\tilde
\lambda'}(y);$$
compare \cite[p.~37]{Macdonald}.
This
identity can be proved by a variant of
Schensted's correspondence, so a bijective proof
of the general equality would be essentially a
generalization of Schensted. Several people have
looked at the problem of a Schensted
generalization corresponding to the case in which
$y_i=0$ when
$i$ is even.
\item[(3)] The analogous results in which $m$ or
$n$ is odd are included in the case in which $m$
and $n$ are both even. For example, if we take
$m=4$ and set $y_3=0$, then the fourth row of a
dimer covering must consist of $n/2$ horizontal
dominoes, which contribute
$\sqrt{x_1x_3\cdots x_{n-1}}$ to the weight, so
we are essentially looking at dimer coverings
with three rows.
\end{enumerate}
\end{problem}

\begin{progress}
A special case of the Robinson--Schensted algorithm given by
Sundquist et al.\ \citeyear{PrSWW}
can be used to get a bijection for a special case of the problem, in
which one sets $y_i = 0$ for all $i$ even, so that we are looking at
dimer coverings (or domino tilings) in which every vertical domino goes
from row $2i+1$ to row $2i+2$ for some $i$. These tilings are not very
interesting because they break up into tilings of 2-by-$n$ rectangles.
But even so, the Robinson--Schensted bijection is nontrivial.
\end{progress}

\section{New Problems}

Let $N(a,b)$ denote the number of matchings
of the $a$-by-$b$ rectangular grid.
Kasteleyn showed that $N(a,b)$ is equal to the square root of
the absolute value of
$$\prod_{j=1}^{a} \prod_{k=1}^{b}
\left(2 \cos \frac{\pi j}{a+1} + 2i \cos \frac{\pi k}{b+1}\right).$$
Some number-theoretic properties of $N(a,b)$ follow from this representation
(see, e.g., \cite{PrCo})
but lack a combinatorial explanation.
The next two problems describe two such facts.

\begin{problem}
\label{prob21}
Give a combinatorial proof of the fact that
$N(a,b)$ divides $N(A,B)$ whenever
$a+1$ divides $A+1$ and $b+1$ divides $B+1$.
\end{problem}

\begin{progress}
Bruce Sagan has given an answer in the ``Fibonacci case'' $a=2$.
A matching of a $2$-by-$(kn-1)$ grid
either splits up as a matching of a $2$-by-$(n-1)$ grid on the left
and a $2$-by-$(kn-n)$ grid on the right
or it splits up as a matching of a $2$-by-$(n-2)$ grid on the left,
a horizontal matching of a $2$-by-$2$ grid in the middle,
and a matching of a $2$-by-$(kn-n-1)$ grid on the right.
Hence
$$N(2,@kn-1) = N(2,@n-1)@N(2,@kn-n) + N(2,@n-2)@N(2,@(k-1)n-1).$$
{}From this formula one can prove that $N(2,n-1)$ divides $N(2,kn-1)$
by induction on $k$.
Volker Strehl has approached the problem in a different way;
his ideas make it seem likely
that a better combinatorial understanding of resultants,
in combination with known interpretations of Chebyshev polynomials,
would be helpful in approaching this problem.
\end{progress}

\begin{problem}
\label{prob22}
Give a combinatorial proof of the fact that
$N(a,2a)$ is always congruent to 1 mod 4.
\end{problem}

\hskip\parindent
(Pachter \citeyear{PrPa} has demonstrated the sort of combinatorial
methods one can use in such problems.)

\medbreak

Even without Kasteleyn's formula,
it is easy to show (e.g., via the transfer-matrix method)
that for any fixed $a$,
the sequence of numbers $N(a,b)$ (with $b$ varying)
satisfies a linear recurrence relation with constant coefficients.
Indeed, consider all $2^a$ different ways
of removing some subset of the $a$ rightmost vertices in the $a$-by-$b$ grid;
this gives us $2^a$ ``mutilated'' versions of the graph.
We can set up recurrences that link
matchings of mutilated graphs of width $b$
with matchings of mutilated graphs of width $b$ and $b-1$,
and standard algebraic methods allow us to turn this system
of joint mutual recurrences of low degree
into a single recurrence of high degree
governing the particular sequence of interest,
which enumerates matchings of unmutilated rectangles.
The recurrence obtained in this way is not, however, best possible,
as one can see even in the simple case $a=2$.

\begin{problem}[Stanley]
Prove or disprove that the minimum degree of a linear recurrence
governing the sequence $N(a,1),N(a,2),N(a,3),\dots$
is $2^{\lfloor (a+1)/2 \rfloor}$.
\end{problem}

\begin{progress}
Observations made by Stanley \citeyear[p.~87]{PrSt2} imply that
the conjecture is true when $a+1$ is an odd prime.
\end{progress}

\hskip\parindent
The idea of mutilating a graph by removing some vertices along its boundary
leads us to the next problem.  It has been observed for small values of $n$
that if one removes equal numbers of black and white vertices from the
boundary of a $2n$-by-$2n$ square grid, the number of matchings
of the mutilated graph is less than the number of matchings of
the original graph.  In fact, it appears to be true that
one can delete \textit{any} subset of the vertices of the square grid and
obtain an induced graph with strictly fewer matchings than the original.

It is worth pointing out that not every graph shares this property with
the square grid.  For instance,
if $G$ is the Aztec diamond graph of order 5 and $G'$ is the graph obtained
from $G$ by deleting the middle vertices along the northwest and northeast
borders, then $G$ has 32768 matchings while $G'$ has 59493.

\begin{problem}
Prove or disprove that every subgraph of the $2n$-by-$2n$ grid graph
has strictly fewer matchings.
\end{problem}

\hskip\parindent
Next we come to a variant on the Aztec dungeon region shown in Figure 4.
Figure 10 shows an ``hexagonal dungeon'' with sides $2,4,4,2,4,4$.
Matt Blum's investigation of these shapes has led him to discover
many patterns; the most striking of these patterns forms the basis
of the next problem.

\begin{figure}
\centerline{\psfig{file=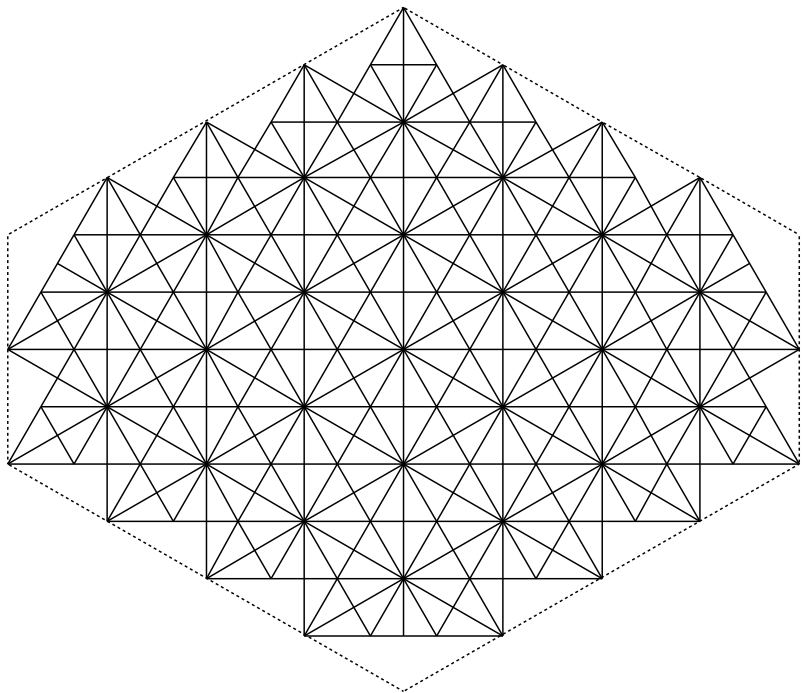,scale=78}}
\caption{An hexagonal dungeon.}
\label{hexdungeon}
\end{figure}

\begin{problem}
Show that the hexagonal dungeon with sides $a,2a,b,a,2a,b$ has exactly
$$13^{2a^2} 14^{\lfloor a^2/2 \rfloor}$$
diform tilings, for all $b \geq 2a$.
\end{problem}

\hskip\parindent
Unmatchable bipartite graphs can sometimes give rise to
interesting quasimatching problems, either by way of $K K^*$
(see Problem 6) or by systematic addition or deletion of vertices or edges.
The former sort of problem simply asks for the determinant
of $K K^*$ (where we may assume that $K$ has more columns than rows).
When the underlying graph has equal numbers of black and white
vertices, this is just the square of the number of matchings,
but when $K$ is a rectangular matrix, $K K^*$ will in general
have a nonzero determinant, even though the graph has no
matchings.

\setcounter{topnumber}{1}

\begin{problem}
Calculate the determinant of $K K^*$ where $K$ is the Kasteleyn--Percus
matrix of the $a,b,c,d,e,f$ honeycomb graph.
\end{problem}

\hskip\parindent
(Note that in this case we can simply take $K$ to be the bipartite
adjacency matrix of the graph.)

Cases of special interest are $a$, $b+1 $, $c $, $a+1 $, $b $, $c+1$
and $a,b,a,b,a,b$ hexagons.  These two cases overlap in the
one-parameter family of $a $, $a+1 $, $a $, $a+1 $, $a $, $a+1$
hexagons.  For instance, in the case of the $3,4,3,4,3,4$ hexagon,
$\det(K K^*)$ is $2^8\cdot 3^3\cdot 7^6$.

\begin{problem}
\label{prob27}
Calculate the determinant of $K K^*$ where $K$ is the Kasteleyn--Percus
matrix of an $m$-by-$n$ Aztec rectangle, or where $K$ is the
Kasteleyn--Percus matrix of the ``fool's diamond'' of order $n$.
(The fool's diamond of order 3 is the following region:
$$\cboard{
  X
 XXX
XXXXX
 XXX
  X
}$$
Fool's diamonds of higher orders are defined in a similar way.)
\end{problem}

\begin{progress}
In the case of Aztec rectangles,
Matt Blum has found general formulas for $\det(KK^*)$ when
$m$ is 1, 2, or 3. For fool's diamonds, we get
$$
\medmuskip 1mu
\displaylines{
1\cr 2\cr 3\cdot 5\cr 2^7\cdot 3\cr 3^2\cdot 5^3\cdot 29\cr
2^9\cdot 3\cdot 5 \cdot 7 \cdot 13^2\cr
7^3\cdot 13^4\cdot 29^2\cr 2^{25}\cdot 3\cdot 7^2\cdot 17^3\cr}
$$
(One might also look at ``fool's rectangles''.)
\end{progress}

\begin{figure}
\centerline{\psfig{file=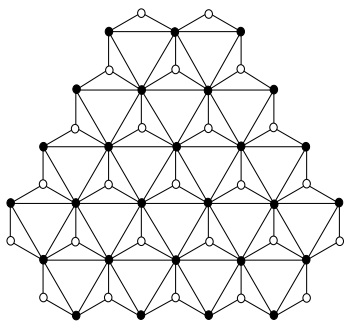,scale=90}}
\caption{A hexagon with extra edges.}
\label{hextwo}
\end{figure}

\hskip\parindent
Another thing one can do with an unmatchable graph is add extra edges.
Even when this ruins the bipartiteness of the graph, there can still
be interesting combinatorics.  For instance, consider the $2,4,2,4,2,4$
hexagon-graph; it has an even number of vertices, but it has a surplus
of black vertices over white vertices.  We therefore introduce
edges between every black vertex and the six nearest black vertices.
(That is, in each hexagon of the honeycomb, we draw a triangle
connecting the three black vertices, as in Figure~\ref{hextwo}.)
Then the graph has $5187 = 3\cdot7\cdot13\cdot19$ matchings.

\begin{problem}
\label{prob28}
Count the matchings of the $a,b,c,d,e,f$
hexagon-graph in which extra edges have been drawn connecting
vertices of the majority color.
\end{problem}

\begin{figure}
\centerline{\psfig{file=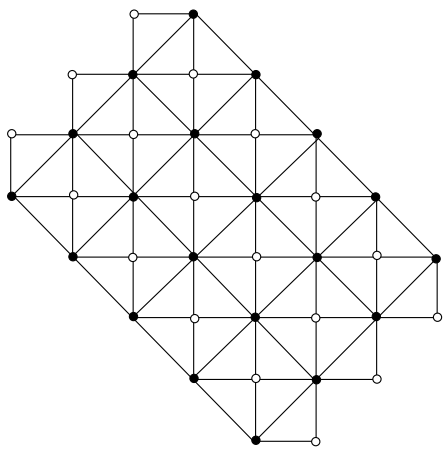,scale=96}}
\caption{An Aztec rectangle with extra edges.}
\label{rectwo}
\end{figure}

\begin{figure}
\centerline{\psfig{file=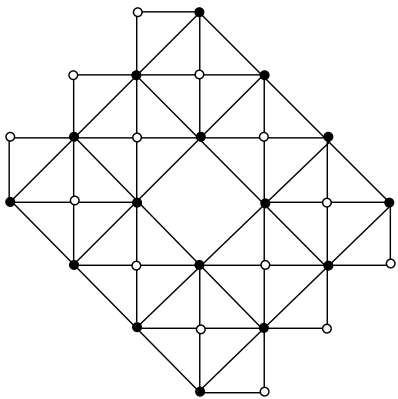,scale=75}}
\caption{A holey Aztec rectangle with extra edges.}
\label{holtwo}
\end{figure}

\hskip\parindent
What works for honeycomb graphs works (or seems to work)
for square-grid graphs as well.  If one adds edges joining each vertex
of majority color to the four nearest like-colored vertices
in the $n$ by $n+2$ Aztec rectangle graph as in Figure~\ref{rectwo},
one gets a graph for which the number of matchings grows like
$2^2\cdot3$, $2^3\cdot3\cdot7$, $2^7\cdot 3\cdot 11$, $2^{17}\cdot 5\cdot 31$,
etc.  If one does the same for the holey $2n-1$ by $2n$ Aztec rectangle
from which the central vertex has been removed, as in Figure~\ref{holtwo},
one gets the numbers
$2^6\cdot 7$, $2^9\cdot 3^2\cdot 13\cdot 17$, $2^{23}\cdot 5^3\cdot 31$, etc.

\begin{problem}
\label{prob29}
Count the matchings of the $a$ by $b$ Aztec
rectangle (with $a+b$ even) in which extra edges have been drawn
connecting vertices of the majority color.  Do the same for the
$2n-1$ by $2n$ holey Aztec rectangle.
\end{problem}

\hskip\parindent
Other examples of nonbipartite graphs for which the number of matchings
has only small prime factors arise when one takes the quotient of a
symmetrical bipartite graph modulo a symmetry-group at least one element
of which interchanges the two colors; Kuperberg \citeyear{PrKu1} gives
some examples of this.
In general, there seem to be fewer product-formula enumerations of
matchings for nonbipartite graphs than for bipartite graphs.
Nevertheless, even in cases where no product formula has been found,
there can be patterns in need of explanation.

Consider the one-parameter family of graphs illustrated in Figure 14
for the case $n=7$ (based on the same nonbipartite infinite graph
\begin{figure}
\centerline{\psfig{file=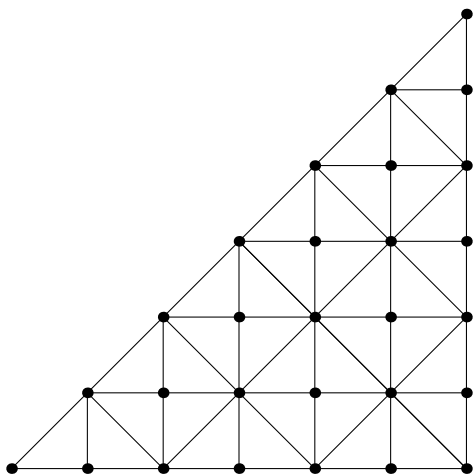,width=2.5in}}
\caption{An isosceles right triangle graph with extra edges.}
\label{diab}
\end{figure}
as Figures 12 and 13).  Such a graph has an even number of vertices
whenever $n$ is congruent to 0 or 3 modulo 4.
Here are the data for the first few cases, courtesy of Matt Blum:
$$
\medmuskip 1mu
\tabskip0pt plus 1fil
\vbox{
\halign to \hsize{\hfil$#$&\hfil#&$#$\hfil\cr
 n&        number of matchings&   $factorization$\cr
 3&                          3&   3\cr
 4&                          6&   2 \cdot 3\cr
 7&                       1065&   3 \cdot 5 \cdot 71\cr
 8&                       6276&   2^2 \cdot 3 \cdot 523\cr
11&                   45949563&   3^2 \cdot 11 \cdot 464137\cr
12&                  807343128&   2^3 \cdot 3^2 \cdot 1109 \cdot 10111\cr
15&            221797080594801&   3^2 \cdot 24644120066089\cr
16&          11812299253803024&   2^4 \cdot 3 \cdot 246089567787563\cr
19&   117066491250943949567763&   3 \cdot 89 \cdot 28289
                                    \cdot 15499002371714201\cr
20& 19100803250397148607852640&   2^5 \cdot 3^2 \cdot 5 \cdot 41 \cdot 367
                                    \cdot 881534305952328473\cr
}}$$

The following problem describes some of Blum's conjectures:

\begin{problem}
\label{prob30}
Show that for the isosceles right triangle graph
with extra edges, the number of matchings is always a multiple of 3.
Furthermore, show that the exact power of 2 dividing the number of
matchings is $2^{n/4}$ when $n$ is 0 modulo 4, and $2^0(=1)$ when
$n$ is 3 modulo 4.
\end{problem}

\hskip\parindent
This property of divisibility by 3 pops up in another problem of
a similar flavor.  Consider the graph shown in Figure 15, which is
just like the one shown in Figure 9, except that half of the triangular
cells have an extra vertex in them, connected to the three nearest
vertices.  (Note also the resemblance to Figure~11.)

\begin{figure}
\centerline{\psfig{file=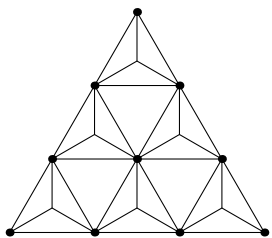,width=1.4in}}
\caption{An equilateral triangle graph with extra vertices and edges.}
\label{newtri}
\end{figure}

\begin{problem}
\label{prob31}
Show that for the equilateral triangle graph
with extra vertices and edges, the number of matchings is always a
multiple of 3.
\end{problem}

\hskip\parindent
(I refrain from making a conjecture about the exponent of 2, though
the data contain patterns suggestive of a general rule.)

\medskip

It may be too soon to try to assemble into one coherent picture all
the diverse phenomena discussed in the preceding 31 problems.  But
I have noticed a gratuitous symmetry that governs many of the exact
formulas, and I will close by pointing it out.  Consider, for example,
the MacMahon--Macdonald product
$$M_n = \prod_{i=0}^{n-1} \prod_{j=0}^{n-1} \prod_{k=0}^{n-1}
\frac{i+j+k+2}{i+j+k+1}$$
that counts matchings of the $n,n,n$ semiregular honeycomb graph.
We find that the ``second quotient'' $M_{n-1} M_{n+1} / M_n^2$
is the rational function
$$\frac{27}{64} \,\frac{(3n-2)(3n-1)^2(3n+1)^2(3n+2)}{(2n-1)^3(2n+1)^3}$$
which is an even function of $n$.

The right hand side in Bo-Yin Yang's theorem (giving the number of diabolo
tilings of a fortress of order $n$) has a power of 5 whose exponent is
$n^2/4$ when $n$ is even and $(n^2-1)/4$ when $n$ is odd; this too is
an even function of $n$.

Domino tilings of Aztec diamonds are enumerated by the formula
$2^{n(n+1)}$.  Here the symmetry is a bit different:
replacing $n$ by $-1-n$ leaves the answer unaffected.

The right hand side of Mihai Ciucu's theorem (giving the number of
diform tilings of an Aztec dungeon of order $n$) has a power of 13
whose exponent is
$(n+1)^2/3$ or $n(n+2)/3$ (according to whether or not $n$ is 2 mod 3).
so that the symmetry corresponds to replacing $n$ by $-2-n$.

There are other instances of this kind that arise,
in which some base is raised to the power of
some quadratic function of $n$;
in each case, the quadratic function
admits a symmetry that preserves the integrality of $n$
(unlike, say, the quadratic function $n(3n+1)/2$,
which as a function from integers to integers
does not possess such a symmetry).

\begin{problem}
\label{prob32}
For many of our formulas,
the ``algebraic'' (right hand) side
is invariant under substitutions that make the
``combinatorial'' (left hand) side
meaningless, insofar as one cannot speak
of graphs with negative numbers of vertices or edges.
Might this invariance nonetheless have some deeper
significance?
\end{problem}

\hskip\parindent
Cohn \citeyear{PrCo} has found
another example of gratuitous symmetry related to tilings.

\section*{Acknowledgements}
This research was conducted with the support of
the National Science Foundation, the National Security Agency,
and the M.I.T. Class of 1922 Career Development fund.
I am deeply indebted to the past and present members
of the Tilings Research Group
for their many forms of assistance:
Pramod Achar, Karen Acquista, Josie Ammer, Federico Ardila,
Rob Blau, Matt Blum, Carl Bosley, Ruth Britto-Pacumio, Constantin Chiscanu,
Henry Cohn, Chris Douglas, Edward Early, Nicholas Eriksson, David
Farris, Lukasz Fidkowski, Marisa Gioioso,
David Gupta, Harald Helfgott, Sharon Hollander, Dan Ionescu,
Sameera Iyengar, Julia Khodor, Neelakantan Krishnaswami,
Eric Kuo, Yvonne Lai, Ching Law, Andrew Menard, Alyce Moy,
Anne-Marie Oreskovich, Ben Raphael,
Vis Taraz, Jordan Weitz, Ben Wieland, Lauren Williams,
David Wilson, Jessica Wong, Jason Woolever, and Laurence Yogman.
I also acknowledge the helpful comments on this manuscript
given by Henry Cohn and Richard Stanley,
and the information provided by Jerry Dias, Michael Fisher and Horst Sachs
concerning the connections between matching theory and the physical sciences.

\bibitemsep=5pt plus 1pt

\def\citeauthoryear #1#2#3{#2\ #3}
\bibliographystyle{msribib}
\bibliography{book37}
\end{document}